\documentclass{article}
\usepackage{eufrak}
\usepackage{euscript}
\usepackage{pictex}
\usepackage{wrapfig}
\usepackage{epsfig}

\def\lnorm{\,\rule[-1mm]{0.6mm}{4mm}\,}

\def\1{\hbox{\upshape1\kern-.15em\vrule height 1.6ex width .3pt
\vrule width .8pt height .25pt\kern.15em}}
\def\der{\,\rule[0mm]{0.1mm}{1.3mm} \rule[0mm]{1.3mm}{0.1mm}
 \hspace{-1.3mm} \rule[1.3mm]{1.3mm}{0.1mm}
 \hspace{-1.3mm} \rule[2.3mm]{1.3mm}{0.1mm}
 \rule[0mm]{0.1mm}{2.3mm} \,}
\def\lev{\langle \,}
\def\des{\,\rangle}

\newtheorem{lemma}{{\sc LEMMA}}[section]
\newtheorem{theorem}{{\sc THEOREM}}[section]
\newtheorem{proposition}{{\sc PROPOSITION}}[section]

\newtheorem{definition}{DEFINITION}[section]
\newtheorem{corollary}{{\sc COROLLARY}}[section]
\newtheorem{example}{{\sc Example}}[section]{\rm}
\newtheorem{Figure}{Figure}[section]

\newcommand {\Proof}  {\mbox{\sc Proof: }}
\newcommand {\QED}{\hfill {\bf \textsf{QED}}}

\newcommand {\Mg} {\overline{M}}
\newcommand {\Md} {\underline{M}}

\def\msx{\mbox{\boldmath{$x$}}}
\def\msxd{\mbox{\boldmath{\scriptsize$x$}}}
\def\msy{\mbox{\boldmath{$y$}}}
\def\msv{\mbox{\boldmath{$v$}}}
\def\msw{\mbox{\boldmath{$w$}}}

\def\msz{\mbox{\boldmath{$z$}}}

\def\msr{\mbox{\boldmath{$r$}}}
\def\msrd{\mbox{\boldmath{\scriptsize$r$}}}
\def\msk{\mbox{\boldmath{$k$}}}
\def\mskd{\mbox{\boldmath{\scriptsize$k$}}}
\def\msm{\mbox{\boldmath{$m$}}}
\def\msmd{\mbox{\boldmath{\scriptsize$m$}}}
\def\msn{\mbox{\boldmath{$n$}}}
\def\msnd{\mbox{\boldmath{\scriptsize$n$}}}
\def\msl{\mbox{\boldmath{$l$}}}
\def\msld{\mbox{\boldmath{\scriptsize$l$}}}
\def\mse{\mbox{\boldmath{$e$}}}
\def\msed{\mbox{\boldmath{\scriptsize$e$}}}
\def\msmu{\mbox{\boldmath{$\mu$}}}

\font\msbm = msbm10
\def\Msbm#1{\hbox{\msbm #1}}

\newcommand {\bbR}{{\Msbm R}}
\newcommand {\bbN}{{\Msbm N}}
\newcommand {\bbZ}{{\Msbm Z}}

\newcommand {\bbJ}{\1}

\newdimen\vskp
\def\AMSclass#1{\vskip\vskp\vbox{\hbox{\vbox{\def\ams{#1}%
	\ifx\ams\empty
	\errmessage{Please put AMS subject classification in \noexpand
	\AMSclass command}
	\else\noindent
	\small\strut{\bf AMS subject classification:\hskip0.5em}\ams
	\strut}}}\fi}

\def\KeyWords#1{\vskip\vskp\vbox{\hbox{\vbox{\def\Kwd{#1}%
	\ifx\Kwd\empty
	\errmessage{Please fill Key words in \noexpand
	\KeyWords command}
	\else
	\newdimen\keywdswidth
	\setbox0=\hbox{\small\bf Key words:\hskip0.5em}
	\keywdswidth=\wd0
	\noindent\hangindent=\keywdswidth
	\hangafter=1
	\small\strut{\bf Key words:\hskip0.5em}\Kwd\strut}}}\fi}

\title{Monotone Numerical Schemes for a Dirichlet Problem for 
Elliptic Operators in Divergence Form}
\author{%
Ned\v zad Limi\'c
	\thanks{%
	Dept. of Mathematics, University of Zagreb, 
	Bijeni\v{c}ka 30, 10002 Zagreb, Croatia,
	e--mail: nlimic@math.hr} 
~and Mladen Rogina	\thanks{%
	Dept. of Mathematics, University of Zagreb, 
	Bijeni\v{c}ka 30, 10002 Zagreb, Croatia,
	e--mail: rogina@math.hr} 
}
\begin{document}
\maketitle 
\begin{abstract}\noindent
We consider a second order differential operator $A(\msx) = -\:\sum_{ i,j=1}^d 
 \partial_i a_{ij}(\msx) \partial_j \:+\: \sum_{j=1}^d \partial_j \big( b_j(\msx) \cdot \big)\:+\: c(\msx)$ on ${\bbR}^d$,
on a bounded domain $D$ with Dirichlet boundary conditions on
$\partial D$, under mild assumptions on the coefficients of the 
diffusion tensor $a_{ij}$. The object is to construct monotone numerical schemes to approximate 
the solution to the problem $A(\msx)\, u(\msx) \: = \: \mu(\msx), \quad \msx \in D$,
where $\mu$ is a positive Radon measure. We start by briefly mentioning questions of 
existence and uniqueness, introducing function spaces needed to prove convergence results.
Then, we define non-standard stencils on grid-knots that lead to extended
discretization schemes by matrices possesing compartmental structure. 
We proceed to discretization of elliptic operators, starting with constant diffusion tensor 
and ending with operators in divergence form. Finally, we discuss $W_2^1$-convergence in detail, and mention
convergence in $C$ and $L_1$ spaces. We conclude by a numerical example illustarting the schemes and
convergence results.
\end{abstract}
\AMSclass{(2000) 35J20, 35J25, 35J15, 65N06, 65N15}
\section{Introduction}\label{sec1}
The object of present analysis are numerical solutions of the elliptic boundary 
value problems in terms of monotone schemes.
It is assumed that the elliptic differential operator has the divergence form
with measurable coefficients satisfying the strict ellipticity condition. The
basic idea of construction of monotone schemes
is presented in \cite{MW} without analysis of convergence of approximate
solutions. Some elaborations of this basic idea can be found in \cite{SMMM, LR2} 
and \cite{LR3} where the convergence is considered in $C$- and $L_1$-spaces. In these
works the stencils of schemes are enclosed by rectangles with vertices at grid-knots.
In the present work we extend the published results by constructing schemes with
stencils stretching far from basic grid-rectangles and so being conceptually 
closer to the original idea of \cite{MW}. The schemes are not derived from finite
difference operators approximating differential operators but rather from a general 
principle which ensures the convergence of approximate solutions. For the case of
a classical elliptic problem this general principle is necessary and sufficient
for the convergence in H\"{o}lder spaces.

In Section \ref{sec2} we describe boundary value problems for linear elliptic 
differential operators which are analyzed here from the numerical point of view.
Problems for bounded domains and Dirichlet boundary conditions are our main
interest. For the sake of completeness, numerical methods for problems on the whole
real space are also studied.

Monotone schemes for linear elliptic problems can be easily constructed by using
discretizations of the corresponding differential operators in terms of matrices 
of positive type. Therefore, in Section \ref{sec3} we start our analysis by
defining matrices of positive type and matrices with the compartmental structure.
The latter ones are more pertinent to construction of discretizations of
elliptic differential operators in divergence form. Monotone schemes give rise to
grid-solutions. By embedding grid-functions into the space of hat functions we
pass from grid-solutions to functions which are called approximate solutions.
The remaining part of this section contains some technical results relating various
norms of grid-solutions and the corresponding norms of approximate solutions.

At the beginning of Section \ref{sec4} we describe discretizations of the operator
$A = -\sum a_{ij}\partial_i \partial_j$ in terms of matrices of positive type by
using forward and backward finite difference operators \cite{LR1}. Then we
describe another method of discretization which is not based on any kind of
finite difference operators. This method becomes a basis for discretizations of
the elliptic operators in divergence form, $A = -\sum \partial_ia_{ij} \partial_j$.
The resulting schemes are called extended schemes. The strong ellipticity of
associated discretized bilinear forms~\cite{Yo} is analyzed and proved for the class
of extended schemes.

The convergence of numerical solutions is analyzed in Section \ref{sec5}. First
we consider the case of $W_2^1$-convergence and later the convergence in H\"{o}lder
spaces and $L_1$-spaces. The convergence in $W_2^1$-space is proved by using
the standard finite element techniques.

Some technical problems arising from a non-continuity of the functions $a_{ij}, i \ne j$,
are discussed in Section \ref{sec6}. It is demonstrated by an example how the extended 
schemes can be applied straightforwardly to problems with general measurable functions
$a_{ij}$.

\section{Preliminaries}\label{sec2}
Elements of ${\bbR}^d$ are denoted by $\msx, \msw, \msr$ {\em etc.}
The Euclidean norm in ${\bbR}^d$ is denoted by $|\cdot|$. For a subset 
$S \subset {\bbR}^d$ the closure is denoted by $\overline{S}$ or $cls(S)$ and boundary
by $\partial S$ or $bnd(S)$. For an open set $D$ the scalar product and norm of $L_2(D)$
are denoted by $(\cdot | \cdot)$ and $\Vert \cdot \Vert_2$, respectively. The norm
of $L_p(D)$ is denoted by $\Vert \cdot \Vert_p$. The Sobolev $W_p^1$-spaces and
$\dot{W}_p^1$-spaces are defined in a
standard way \cite{Ma,St}. In addition we need the H\"{o}lder spaces 
$C^{(\alpha)}(\overline{D}), \alpha \in (0,1)$, defined as in \cite{Ma,St}. 
The convex set of positive Radon measures $\mu$ on $\mathfrak{B}(D)$,
$\int d\mu = 1$ (tight probability measures) is denoted by
${\cal P}(D)$. Then $\langle v \,|\, \mu \rangle = \int_D v(\msx) \mu(d \msx)$ 
is well defined for $v \in \dot{W}_\infty^1(D)$. We say that a sequence of
$\mu_n \in {\cal P}(D)$ converges weakly to $\mu \in {\cal P}(D)$ if $\lim_n \lev v |
\mu_n \des = \lev v | \mu \des$ for each $v \in \dot{C}(\overline{D})$.
An open subset of ${\bbR}^d$ with Lipshitz boundary \cite{St} is called 
a domain with Lipshitz boundary. Only such domains are considered here. The
following result is used. Let $\dot{C}^{(1)}(\overline{D})$ be the linear space of 
continuous functions, having continuous partial derivatives and zero values on
$\partial D$. Then, for a domain $D$ with Lipshitz boundary, the space
$\dot{C}^{(1)}(\overline{D})$ is dense in $\dot{W}_p^1(D), 1 \leq p < \infty$ \cite{Ma}.

Here we study elliptic operators on ${\bbR}^d$ in divergence form defined by
\begin{equation}\label{exp2.2}
 A(\msx) = -\:\sum_{ i,j=1}^d  \partial_i a_{ij}(\msx) \partial_j
 \:+\: \sum_{j=1}^d \partial_j \big( b_j(\msx) \cdot \big)\:+\: c(\msx),
\end{equation}
for which the coefficients are constrained as follows. The functions  
$a_{ij}= a_{ji}$, $b_{i}, i,j = 1,2,\dots, d$ and $c$ are measurable
on ${\bbR}^d$, $c \geq 0$ and $a_{ij}(\msx)$ converge to constant values as $|\msx|$
increases. Apart from this, the elliptic operator must be strictly elliptic,
meaning that there are positive numbers $\underline{M}, \overline{M}, \:
0 < \underline{M} \leq \overline{M}$, such that the double inequality
\begin{equation}\label{exp2.1}
  \underline{M} \,|\msx|^2  \,\leq\,\sum_{i,j=1}^{d}  a_{ij}(\msx)
z_i \bar{z}_{j}
 \leq \,\overline{M}\, |\msx|^2 ,  \quad  \msx \in  {\bbR}^d,
\end{equation}
holds. The main part of $A$ is denoted by $A_0$.

For sufficiently large and positive $\lambda$ there exist kernels
$(\msx, \msy) \mapsto K_A(\lambda,\msx, \msy)$ on $D \times D$, such that
$(\msx, \msy) \mapsto |\msx-\msy|^{d-2}K_A(\lambda,\msx, \msy)$
are uniformly bounded on $D \times D$ and H\"{o}lder continuous in $\msx$
on the set $D \setminus \msy$. The integral operator  with kernel 
$K_A(\lambda,\cdot, \cdot)$ is denoted by $R(\lambda,A)$.
It maps measurable and bounded 
functions with compact supports into continuous functions on $D$ and
$(\lambda I+A)R(\lambda,A)=I$ in $C(\overline{D})$. Let
$\mathfrak{L}(D)$ be the linear space of functions on $D$ and
$\mathfrak{D}(A_D) = R(\lambda, A) \mathfrak{L}(D)$. Then
$A(\msx): \mathfrak{D}(A_D) \mapsto \mathfrak{L}(D)$
is denoted by $A_D$. We say that $A_D$ is defined by $A(\msx)$ 
on $D$ with Dirichlet boundary conditions on $\partial D$.
If $\mathfrak{L}(D)=L_p(D)$, $p < \infty$, the operator $A_D$ is maximally
closed from $\mathfrak{D}(A_D) \subset L_p(D)$ into $L_p(D)$.

Let us define a real bilinear form on $W_q^1(D) \times W_p^1(D)$, $1/p+1/q = 1$, by:
\begin{equation}\label{ex2.3}\begin{array}{c}\displaystyle
 a(v,u) \ = \ \sum_{i,j =1}^d\: \int_D\: a_{ij}( \msx)\:
 \partial_i v(\msx)\:  \partial_j u(\msx) \:d\msx \\
 \displaystyle 
 - \sum_{i=1}^d \: \int_D\: b_i(\msx) \:\partial_i v(\msx)\, 
 u(\msx) \:d \msx +  \int_D\: c(\msx)\, v(\msx)\,  u(\msx)\,d\msx .  \end{array}
\end{equation}
For a domain $D$ with Lipshitz boundary $\partial D$ and for each pair
$v \in \dot{W}_q^1(D)$, $u \in \dot{W}_p^1(D) \cap \{A u \,\in \, (\dot{W}_q^1(D))^{\dag}\} $, $1<p<\infty$,
there holds the Green's formula:
\[ a(v,u) \ = \ \lev v\,|\,A\, u \des . \]
The Green's formula is also valid for each pair
$v \in \dot{W}_{\infty}^1(D)$, $u \in \dot{W}_1^1(D)\cap \{ Au \in {\cal P}(D)\}$.

The boundary value problem, to be studied in this work, is defined by
\begin{equation}\label{exp2.3}\begin{array}{c}
 A(\msx)\, u(\msx) \: = \: \mu(\msx), \quad \msx \in D,\\
 u\,\big | \,\partial D \ = \ 0, \end{array}
\end{equation}
where $D$ is a domain with Lipshitz boundary and the nonhomogeneous term $\mu$ is either
a positive Radon measure or $\mu \in \dot{W}_2^{-1}(D)$. In the case of $D = {\bbR}^d$ we understand
that the boundary condition in (\ref{exp2.3}) is omitted.

The variational formulation of (\ref{exp2.3}) for a solution
$u \in \dot{W}_p^1(D)$, $p=1,2$, has the following form:
\begin{equation}\label{ex2.5}
 a(v,u) \ = \ \langle \,v \,|\,\mu \,\rangle, \quad \textrm{for~any} \
 v \in \dot{W}_q^1(D).
\end{equation}
Solutions of (\ref{exp2.3}) and (\ref{ex2.5}) are called strong and weak
solutions, respectively. 

The following result is valid \cite{LR2}.

\begin{theorem}\label{th2.1} Let $D$ be a bounded domain with Lipshitz
boundary. For each $p \in [1, d/(d-1))$ there exists
a unique weak solution $u$ of (\ref{ex2.5}) belonging to the class
$\dot{W}_p^1(D)$ and possessing the following properties:
\begin{description}\itemsep 0.cm
 \item{(i)} There exists a positive number $c$ depending on
$\underline{M}, \overline{M}, p, D$, such that the following inequality
is valid:
\[ \Vert \,u\, \Vert_{p,1} \ < \ c.  \]
 \item{(ii)} If $\{ \mu_n : n \in {\bbN}\} \subset {\cal P}(D)$
converges weakly to a $\mu \in {\cal P}(D)$, then the corresponding sequence
of weak
solutions $\{ u_n : n \in {\bbN}\} \subset \dot{W}_p^1(D)$, $u_n =
A_D^{-1}\mu_n$,
converges strongly in $L_p(D)$ to $u = A_D^{-1} \mu$.
\end{description}
\end{theorem}

\section{Grid-functions and $W_2^1$-spaces}\label{sec3}
Let the orthogonal coordinate system in ${\bbR}^d$ be determined by unit vectors
$\mse_i$, and let us, for each $n \in {\bbN}$, define a numerical grid $G_n$ on 
${\bbR}^d$ by vectors $\msx = \sum_{l=1}^d h(n)\, k_l \mse_l$, where the grid-step
$h(n)$ is determined by $h(n) = 2^{-n}$. If not necessary the grid-steps 
$h(n)$ are shortly denoted by $h$. For the sake of simple notation the grid
step $h$ is assumed to be the same for all coordinate directions. The obtained results
of convergence are valid as well for grids with grid-steps depending on direction
$\mse_i$. To each $\msv \in G_n$ there corresponds a
grid-cube $C(h, \msv) = \prod_1^d \,[v_j,v_j+h)$, where $v_j$ are
coordinates of $\msv$. Cubes $C(h,\msv)$ define
a decomposition of ${\bbR}^d$ into disjoint sets. We say that 
the sets $G_n(D) = G_n \cap D$ are {\em discretizations} of $D$.
The grids $G_n$ are homogeneous with respect to
translations in the direction of coordinate axes, {\em i.e.} for $\msx \in G_n$ 
and $\mbox{\boldmath{$t$}}=h p_i\,\mse_i, p_i \in {\bbZ}$ we have
$\msx + \mbox{\boldmath{$t$}} \in G_n$. There exist subsets of $G_n$ which are 
also homogeneous in the defined sense. Let $\msr_0 \in G_n$ and 
$\msr = (r_1, r_2, \ldots, r_d) \in {\bbN}^d$ be fixed. The set
\begin{equation}\label{exn2.2}
 G_n(\msr_0, \msr) \ = \  \{ \msr_0 \:+\: h\,\sum_{l=1}^d \:k_l\,r_l \,
 \mse_l \ : \ k_l \in {\bbZ}\} 
\end{equation}
is a subsets of $G_n$. It is homogeneous with respect to translations for
$hk_ir_i\mse_i,\: k_i \in {\bbZ}$, grid-knots. There are 
$vol(R)=\prod_{i=1}^d r_i$ disjoint subgrids (\ref{exn2.2}) making a partition of
$G_n$. Each is denoted $G_n(R)$, where $R$ stands shortly for $2d$ parameters
$\msr_0 \in G_n, \msr \in {\bbN}^d$. The index set of $G_n$ is denoted 
by $I_n$. Similarly, $I_n(R)$ is the index set of $G_n(R)$.

The discretization of a function $u \in C({\bbR}^d)$ on $G_n$ is denoted by 
${\bf u}_n$ and defined by values at grid-knots, $\big ({\bf u}_n\big )_{\msmd} \:=\:
u(\msx_{\msmd})$ where $\msx_{\msmd} = (m_1h, m_2h, \ldots,m_dh) \in G_n$, 
and $\msm = (m_1, m_2, \ldots,m_d)$ is a multi-index. The function ${\bf u}_n$ is
usually called a grid function. We denote the linear spaces 
of discretizations by $l(G_n)$ or $l(G_n(D))$. Elements of $l(G_n)$ are also
called columns. The corresponding $L_p$-spaces are denoted by $l_p(G_n)$ or 
$l_p(G_n(D))$, and their norms by $\lnorm \cdot \lnorm_p$. The duality pairing
of ${\bf v} \in l_q(G_n)$ and ${\bf u} \in l_p(G_n)$ is denoted by $\lev
{\bf v} |{\bf u} \des$. The scalar product in $l_2(G_n)$ is denoted by
$\lev \cdot |\cdot \des$ and sometimes by $( \cdot |\cdot )$. 
The norm of $l_p(G_n(R))$ is denoted by $\lnorm \cdot \lnorm_{Rp}$. 
For $p \in [1,\infty)$ we have
\[ \lnorm {\bf u} \lnorm_{Rp} \ = \ \Big[ \,vol(R)\,\sum_{\mskd \in I_n(R)}\:
 |u_{\mskd}|^p\, \Big]^{1/p}, \]
while $\lnorm {\bf u} \lnorm_{R\infty} = \sup\{|u_{\mskd}|: \msk \in I_n(R)\}$.
A matrix $A_n$ with indices 
corresponding to grid-knots of $G_n$ (or $G_n(D)$) is said to be defined on 
$G_n$ (or $G_n(D)$). We use the notation $A_n \geq 0$
if all entries of $A_n$ are non negative.

The {\em shift operator} $Z(\msx), \msx \in {\bbR}^d$, acting on functions 
$f:{\bbR}^d \mapsto {\bbR}$, is defined by $\big(Z(\msx)f\big)(\msx) = 
f(\msx+\msz)$. Similarly we define the discretized shift operator by 
$\big ( Z_n(r,i) {\bf u}_n \big)_{\msmd} = ({\bf u}_n)_{\msnd}$, where $\msn = 
\msm + rh\mse_i$.

{\bf Discretization of differential operators}. 
A function $u \in C^{(1)}({\bbR}^d)$ has continuous partial derivatives
$\partial_i u, i = 1,2,\ldots.d$. With respect to a grid step $h$, the partial
derivatives are discretized by forward/backward finite difference operators in the 
usual way:
\begin{equation}\label{ex2.7}\begin{array}{c}
 \der_i(t) u(\msx) \ = \ \frac{1}{t} \big( 
 u(\msx + t \mse_i) \:-\:
 u(\msx) \big ),\\
 \widehat{\der}_i(t) u(\msx) \ = \ \frac{1}{t} \big( 
 u(\msx) \:-\:
 u(\msx - t\mse_i) \big ). \end{array}, \quad \msx \in {\bbR}^d, \ t \ne 0.
\end{equation}
Discretizations of the functions $\partial_iu$ on $G_n$, denoted by
$U_i(r){\bf u}_n, V_i(r){\bf u}_n$, are defined by:
\[ \big (U_i(r)\,{\bf u}_n\big )_{\msmd} \:=\: \der_i\big(rh \big)
 \:u(\msx_{\msmd}), \quad  \big (V_i(r)\,{\bf u}_n\big )_{\msmd} \:=\:
 \widehat{\der}_i \big(rh \big) \:u(\msx_{\msmd}).\]
Then 
\[ \begin{array}{l} U_i(r) \ = \ (rh)^{-1}(Z_n(r,i) \,-\, I),\\
 V_i(r) \ = \  (rh)^{-1} \big(I - Z_n(-r,i) \big) \:= \:
 U_i(-r) \:=\:\,-\, U_i(r)^T.
\end{array}\]
Therefore we have $U_i(-r) = U_i(r)\,Z_n(-r,i) = Z_n(-r,i)\,U_i(r)$,
and similarly for $V_i(r)$. 

In accordance with the previous terminology, we say that $\partial_i, \sum_{ij}
\partial_i a_{ij} \partial_j$ {\em etc.} are differential operators on ${\bbR}^d$ or $D$.
We say that their discretizations are defined on $G_n$ or $G_n(D)$. In
particular, discretizations of the differential operator (\ref{exp2.2}) are
denoted by $A_n$. Naturally, matrices $A_n$ are the main object in this work.
We intend to analyze a class of discretizations $A_n$ of 
(\ref{exp2.2}) with compartmental structure. 

\begin{definition}[Compartmental structure]\label{defn1.1}
A matrix $A = \{a_{ij}\}_{II}$ is said to have the compartmental
structure if it has positive diagonal entries, non-positive off-diagonal 
entries and positive or zero column sums:
\[ a_{ii} \:\geq \: 0, \qquad a_{ij} \:\leq \: 0 {\rm~~for~~}i \,\ne \,
 j, \qquad a_j \:=\: \sum_{i \in I} \:a_{ij} \:\geq \:0.\]
A matrix with the compartmental structure with zero column sums 
is called conservative. If $A$ is compartmental then $B = A^T$ is called
a matrix of positive type. A matrix of positive type $B$ is called 
conservative if $A = B^T$ is (compartmental) conservative matrix.
\end{definition}

Hence, we consider a class of discretizations $A_n$ with the
following properties:
\begin{equation}\label{exn11.3}\begin{array}{lll}
 (A_n)_{\mskd\mskd} \ > \ 0, & (A_n)_{\mskd\msld} \ \leq 0,&
 \msl \ne \msk,  \\
 \sum_{\mskd}\: (A_n)_{\mskd \msld} \ \geq \ 0. &&\end{array}
\end{equation}
Additionally, there must exist a positive number $\sigma^2$ such that
\[ \sigma^2 \ = \ \sup_{n} \:h(n)^2\: \sup\,\{ (A_n)_{\mskd\mskd} \ : \
 \msk \in I_n\}.\]
Then $Q_n = \sigma^2I -h^2A_n \geq 0$ and $\lnorm Q_n \lnorm_1 \leq \sigma^2$.
This fact enables a decomposition of $A_n$ as $A_n = h^{-2}\big( \sigma^2 I
- Q_n\big)$. Hence, the resolvents
\[ R(\lambda,A_n) \ = \ \frac{h^2}{\sigma^2}\: \frac{1}{1 + \lambda^2/\sigma^2}\:
 \sum_{k=0}^\infty \:\frac{1}{(1 + \lambda^2/\sigma^2)^k}\:\, Q_n^k,  \]
are non-negative matrices on $G_n$, $\lnorm R(\lambda,A_n)\lnorm_1 \leq 1/
\lambda$.

Let $I' \subset I$. Then $A' = \{a_{ij}\}_{I'I'}$ is called 
a diagonal submatrix of $A=\{a_{ij}\}_{II}$. A non-negative matrix
$A = \{a_{ij}\}_{II}$ is called irreducible if for each finite index set $I' \subset I$
the corresponding diagonal submatrix is irreducible.
A non-negative matrix which is not irreducible is called reducible.
Thus a non-negative matrix $A$ is irreducible if for each finite subsets $I' 
\subset I$ there exists $m \in {\bbN}$ such that the entries $p_{ij}$ of $P =
A^m$ are positive for $i,j \in I'$.

\begin{lemma}\label{lemn11.1} Let $A_n$ possess either the compartmental structure or
be of positive type. Then the matrix $R(\lambda,A_n), \lambda > 0$ on $G_n$ 
is positive iff it is irreducible.
\end{lemma}

Let the matrices $A_n$ on $G_n$ be discretizations of the differential operator
(\ref{exp2.2}) and let there be defined linear systems:
\begin{equation}\label{ex2.10}
 A_n {\bf u}_n \ = \ \msmu_n,
\end{equation}
where ${\bf u}_n, \msmu_n \in l(G_n(D))$. We say that the system (\ref{ex2.10})
numerically approximates the boundary value problem (\ref{exp2.3}). The columns
${\bf u}_n$ are called grid-solutions. Obviously, $\msmu_n$ are discretizations of
$\mu$. If the matrix $A_n$ is compartmental
or of positive type we have $A_n = D_n - B_n$, where the diagonal matrix $D_n$ has
positive entries and non-negative matrix $B_n$ has zero diagonal entries. Finite difference
equations (\ref{ex2.10}) can be rewritten as ${\bf u}_n = D_n^{-1}B_n{\bf u}_n +
D_n^{-1}\msmu_n$ or componentwise as
\[ u_{\mskd} \ = \ \sum_{\msld \ne \mskd}\: \big(D_n^{-1}B_n\big)_{\mskd \msld}\,
 u_{\msld} \:+\: d_{\mskd \mskd}^{-1}\,\msmu_{\mskd}.\]
Hence, if $A_n$ is compartmental or of positive type, then the obtained numerical 
scheme is monotone~\cite{BS}.

For a matrix $A_n$ on $G_n$ we define numerical neighborhoods
\[  {\cal N}(\msx) \ = \ \{ \msy \in  G_n \::\: \msx = h\msk,\:
 \msy = h\msl, \ (A_n)_{ \mskd \msld} \ne 0\}. \]
Obviously, numerical neighborhoods of a system matrix and stencils of the
corresponding finite difference schemes are mutually related.

\subsection*{Imbedding of grid-functions into $W_2^1$-space}
Let us define the quadratic functional on $l(G_n)$ by $q({\bf u}) =
\sum_i^d \lnorm U_i{\bf u}\lnorm_2^2$ and
$q_R({\bf u}) = vol(R)\sum_i^d \lnorm U_i(r_i){\bf u}\lnorm_{R2}^2$ on $l(G_n(R))$.
It is understood that $q_R = q$ for $G_n(R) = G_n$. There exist symmetric matrices $Q_n$ 
on $G_n$ such that $q_R({\bf u}_n) = \lev {\bf u}|Q_n {\bf u}\des$.
Discrete analogs of $W_2^1$-spaces are spaces of those 
${\bf u}_n \in l(G_n(R))$ for which the norm $\lnorm \cdot \lnorm_{R2,1}$:
\begin{equation}\label{exn2.12}
 \lnorm {\bf u} \lnorm_{R2, 1}^2 \ = \ \lnorm {\bf u} \lnorm_{R2}^2
 \ + \ q_R({\bf u}),
\end{equation}
is finite. By convention $\lnorm \cdot \lnorm_{2,1} = \lnorm \cdot \lnorm_{R2,1}$
for $r_i = 1$. The subspace of grid-functions ${\bf u} \in w_2^1(G_n(R))$ for 
which ${\bf u}_n = {\bbJ}_{G_n(D)}{\bf u}_n$ is denoted by $w_2^1(G_n(R,D))$. 
Hence, $w_2^1(G_n(R,D))$ for $\msr = {\bf 1}$ is denoted by $w_2^1(G_n(D))$. The
restriction of $q_R$ on $G_n(R,D)$ is represented as $q_R({\bf u}) =
\lev {\bf u}|Q_n(D){\bf u}\des$, where $Q_n(D)$ is a symmetric matrix on
$G_n(R,D)$. By using the fact that the negative Laplacean
on $\dot{C}^{(2)}(\overline{D})$ and its discretizations $-\sum_i V_iU_i$ in
$l_2(G_n(D))$ have positive minimal eigenvalues, we can derive the following
result.

\begin{lemma}\label{lem3.1} Let $D$ be bounded and $Q_n(D)$ be irreducible on 
$G_n(R,D)$. Then the norms (\ref{exn2.12}) and $q_R(\cdot)^{1/2}$ 
are equivalent in $w_2^1(G_n(R,D))$, 
\[ q_R(\cdot)^{1/2} \ \geq \ \beta\,\lnorm \cdot \lnorm_{R2,1},\]
where $\beta$ is independent of $n$.
\end{lemma}

Let us consider a norm $\lnorm \cdot \lnorm_{R2,1}$ on $l(G_n(R))$ defined by
(\ref{exn2.12}). Any such norm is a semi-norm on $l(G_n)$. Our 
object of interest are quadratic functionals:
\begin{equation}\label{exn2.13}\begin{array}{lll} \displaystyle
 \lnorm {\bf u} \lnorm_{avg,2,1}^2 &=& \frac{1}{vol(R)}\,\sum_{R}\:
 \lnorm {\bf u} \lnorm_{R2,1}^2 \\ \displaystyle
 &=&  \lnorm {\bf u} \lnorm_{2}^2\:+\:\frac{1}{vol(R)}
 \sum_{i=1}^d\,\lnorm U_i(r_i){\bf u} \lnorm_{R2}^2 \ \leq \  \lnorm {\bf u} \lnorm_{R2,1}^2
.\end{array}  \ {\bf u} \in l(G_n).
\end{equation}
Then $\lnorm \cdot \lnorm_{avg,2,1}$ is a norm on $l(G_n)$. Unfortunately, it is not
equivalent to $\lnorm \cdot \lnorm_{2,1}$ uniformly with respect to $n$.

An element (column) ${\bf u}_n \in l(G_n)$ can be associated to a continuous function
on ${\bbR}^d$ in various ways. Here is utilized a mapping $l(G_n) \mapsto C({\bbR}^d)$
which is defined in terms of hat functions. Let $\chi$ be the canonical hat function
on ${\bbR}$, centered at the origin and having the support $[-1,1]$. Then
$z \mapsto \phi(h,x,z) = \chi(h^{-1}(z-hx))$ is
the hat function on ${\bbR}$, centered
at $x \in {\bbR}$ with support $[x-h,x+h]$. The functions
$\msz  \: \mapsto \: \phi_{\mskd}(\msz) \:=\:\prod_{i=1}^d \:\phi(h,x_i,z_i),
x_i = hk_i$, define $d$-dimensional hat functions with supports 
$S({\bf 1},\msx) = \prod_i [x_i-h, x_i+h]$. The functions $\phi_{\mskd}(\cdot)
\in G_n$, span a linear space, denoted by $E_n({\bbR}^d)$.
Let ${\bf u}_n \in l(G_n)$ have the entries $u_{n \mskd} = ({\bf u}_n)_{\mskd}$.
Then the function $u(n) =\sum_{\mskd \in {\bbZ}^d} u_{n  \mskd} \phi_{\mskd}$
belongs to $E_n({\bbR}^d)$ and defines imbedding of grid-functions into the
space of continuous functions. We denote the corresponding mapping by
$\Phi_n : l(G_n) \mapsto E_n({\bbR}^d)$. Obviously there exists
$\Phi_n^{-1} : E_n({\bbR}^d)  \mapsto l(G_n)$ and the spaces 
$l(G_n)$ and $E_n({\bbR}^d)$ are isomorphic with respect to the pair
of mappings $\Phi_n, \Phi_n^{-1}$. It is clear that $E_{n}({\bbR}^d)
\subset E_{n+1}({\bbR}^d)$ and the space of functions $\cup_n E_n({\bbR}^d)$
is dense in $L_p({\bbR}^d), p \in [1,\infty)$, as well as in
$\dot{C}({\bbR}^d)$. Let us mention that $\sum_{\mskd}
\phi_{\mskd} = 1$ on ${\bbR}^d$.

Now we consider another collection of basis functions. To each $\msx = h\msk \in
G_n(R)$ there is associated a $d$-dimensional hat function
\[ \psi_{\mskd}(\msx) \ = \ \prod_{i=1}^d\: \chi \left(\frac{x_i-hk_i}{hr_i}\right ),\]
obviously, with the support $S(\msr,\msx) = \prod_i [x_i-r_ih, x_i+r_ih]$. They
span a linear space denoted by $E_n(R,{\bbR}^d)$. Again
we have $\sum_{\mskd}\psi_{\mskd} = 1$ on ${\bbR}^d$.
The mappings $\Phi_n, \Phi_n^{-1}$ cannot be applied to elements
of $l(G_n(R))$ and $E_n(R,{\bbR}^d)$, respectively. Therefore we define
restrictions $\Phi_n(R): l(G_n(R)) \rightarrow E_n(R,{\bbR}^d)$ and
$\Phi_n^{-1}(R)$ by the following expression:
\begin{equation}\label{exn3.5}
 u(n) \ = \ \Phi_n(R)\,{\bf u}_n \ = \ \sum_{\mskd}\:
 \big({\bf u}_n\big)_{\mskd}\:\psi_{\mskd}.
\end{equation}
If we have to emphasize that $u(n)$ is related to a particular set of parameters
$R$, then we use an extended denotation $u(R,n)$.

The following results are necessary in our proof of the consistency in Section \ref{sec5}:

\begin{lemma}\label{lem3.2} Let for each $n \in {\bbN}$ there exist 
columns ${\bf u}_n,{\bf v}_n \in l(G_n(R))$ such that 
$h^d \lnorm {\bf u}_n \lnorm_{R2}^2 \leq c, h^d \lnorm {\bf v}_n
\lnorm_{R2}^2 \leq c$, where $c > 0$ is independent of $n$.
Let $u(n) = \Phi_n(R) {\bf u}_n$ and $v(n) = \Phi_n(R) {\bf v}_n$. If 
\[ \lim_n \:h^d\: \sup_{|\mskd-\msld| =|\msrd|} 
 \lnorm (Z_n(\msk-\msl)  - I) {\bf u}_n \lnorm_{R2}^2 \ = \ 0, \]
then
\[ \lim_n \Big (\,(v(n) | u(n)) \ - \ h^d\,({\bf v}_n\,|\,{\bf u}_n)_R
 \,\Big ) \ = \ 0.\]
\end{lemma}

It is easy to check $\Vert u(n) \Vert_2^2 \leq h^d\lnorm {\bf u}_n\lnorm_2^2$.
There holds a similar inequality for any pair $\partial_iu, U_i{\bf u}_n$. First
one verifies 
\[ (\partial_iv(n)\,|\,\partial_iu(n)) \ = \ (\der_i(r_ih)v(n)\,|\,\der_i(hr_i)(u(n))
 \ = \ (U_i(r_i)\,{\bf v}_n \,|\,U_i(r_i)\,{\bf u}_n),\]
so that the following Theorem holds.

\begin{theorem}\label{thn3.1} Let sequences of functions $v(n), u(n), 
n \in {\bbN}$, be defined by (\ref{exn3.5}). Then 
\[\begin{array}{lc} (i)&
 \Big| \sum_{i=1}^d \: \big( \,\partial_i v(n)\,|\,\partial_i u(n)\,\big)
 \Big | \ \leq \  h^d\,q_R({\bf v}_n)^{1/2}\,q_R({\bf u}_n)^{1/2}.\\
 (ii)& 
 \Big|\big( \,\partial_i v(n)\,|\,\partial_j u(n)\,\big) \,-\,h^d\,
 \sum_{\mskd} \:\big(U_i(r_i){\bf v}\big)_{\mskd}\,\big(U_j(r_j){\bf u}\big)_{\mskd} 
 \Big|\\ 
 &\leq \ h^d \, \min\:\left \{ \begin{array}{l} \displaystyle
 \lnorm U_i(r_i){\bf v}\lnorm_2\,\sup \:\big\{\lnorm (Z(w,j)-I)\,
 U_j(r_j){\bf u}\lnorm_2 \::\: |w| \leq r_jh \big\} 
 \\ \displaystyle
 \lnorm U_j(r_j){\bf u}\lnorm_2\,\sup \:\big\{\lnorm (Z(w,i)-I)\,
 U_i(r_i){\bf v}\lnorm_2 \::\: |w| \leq r_ih \big\} 
 \end{array} \right.  \end{array}\]
\end{theorem}

Now we can get the following useful result involving the norm $\Vert u(n)\Vert_{2,1}$
and its averaged value:
\begin{equation}\label{ex3.4}
 \Vert u(n) \Vert_{avg,2,1}^2 \ = \  \frac{1}{vol(R)}\,\sum_R\:\Vert u(R,n) 
 \Vert_{2,1}^2.
\end{equation}

\begin{corollary}\label{corn3.1} Let $u(n) = \Phi_n {\bf u}_n$. There exists $\sigma^2
\in (0,1)$, independent of $n$, such that
\[\begin{array}{lcccl}\displaystyle
 (1-\sigma^2)\,h^d\,\lnorm {\bf u}_n \lnorm_{R2,1}^2 &\leq&
 \Vert u(R,n) \Vert_{2,1}^2 &\leq& h^d \, \lnorm {\bf u}_n \lnorm_{R2,1}^2,\\ \displaystyle
 (1-\sigma^2)\,h^d\,\lnorm {\bf u}_n \lnorm_{avg,2,1}^2 &\leq&
 \Vert u(n) \Vert_{avg,2,1}^2 &\leq& h^d \, \lnorm {\bf u}_n \lnorm_{avg,2,1}^2.
 \end{array}\]
\end{corollary}

An element $u \in W_2^1({\bbR}^d)$ does not belong necessary to $E_n(R,{\bbR}^d)$.
In order to approximate $u$ with elements of $E_n(R,{\bbR}^d)$ we define:
\begin{equation}\label{exn3.3}
 \hat{u}(n) \ = \ \sum_{\mskd \in I_n(R)}\: \Vert \psi_{\mskd} \Vert_1^{-1}
 \:(\psi_{\mskd}|u)\,\psi_{\mskd}.
\end{equation}
The numbers $\Vert \psi_{\mskd} \Vert_1^{-1}(\psi_{\mskd}|u)$ are called 
Fourier coefficients of $u$. 

The basic result for our proof of convergence of approximate solutions is
formulated by using the quantity
$\Gamma_p(\msw,u)$ defined by:
\[ \Gamma_p(\msw,u) \ = \  \Vert \, (Z(\msw) \,-\,I) \,u\, \Vert_p .\]
The kernels
\[ \omega_n(\msx,\msy) \ = \ \sum_{\mskd} \: \frac{1}{\Vert \psi_{\mskd}\Vert_1}\:
 \psi_{\mskd}(\msx)\,\psi_{\mskd}(\msy)\]
define an integral operator which is denoted by $K_n$. Actually, the kernels
$\omega_n$ define a $\delta$-sequence of functions on ${\bbR}^d \times {\bbR}^d$
and $K_n$ converge strongly in $L_p$-spaces to unity:

\begin{corollary}\label{corn3.0} Let $p \in [1,\infty]$. Then
\begin{description}\itemsep 0.cm
 \item{(i)} $\Vert K_n \Vert_p \:\leq \:1$.
 \item{(ii)} There is a positive number $\kappa(R)$, independent of $n$, such
that $\Vert (I-K_n)u \Vert_p \:\leq \: \kappa(R)
 \sup\{\Gamma_p(\msw,u)\,:\,|w_i| \leq \:hr_i\}$.
 \item{(iii)} The operator $K_n \in L(L_2({\bbR}^d),L_2({\bbR}^d))$ has the
spectrum equal $[0,1]$.
\end{description}
\end{corollary}

If necessary, we write $K_n(R)$ instead of $K_n$ to emphasize that the integral
operator is defined with a basis $\psi_{\mskd}$ in $E_n(R,{\bbR}^d)$ for a particular
value of parameter $R$.

\begin{theorem}\label{Thn3.2} Let $v,u \in W_2^1({\bbR}^d)$ and 
$\hat{u}(n),\hat{v}(n)$ be defined by (\ref{exn3.3}). Then
\[\begin{array}{c}
 \Big |(\hat{v}(n) |\hat{u}(n)) \:-\: (v|u)\Big | \ \leq \ 
 c(R)\,\min \:\left \{\begin{array}{l}
 \Vert \,u\,\Vert_2 \:  \sup_{|\msw| = h}\:\Gamma_2(\msw,v), \\
 \Vert \,v\,\Vert_2 \:  \sup_{|\msw| = h}\:\Gamma_2(\msw,u), \\
 \end{array} \right . \\  \displaystyle
 \Big |(\partial \hat{v}(n)\,|\,\partial \hat{u}(n)) \:-\:
 (\,\partial v\,|\, \partial u\,)\Big | 
 \leq \  c(R)\,\min\:\left \{ \begin{array}{l} \displaystyle
 \Vert \partial v \Vert_2\,\left [  \Vert \partial u - 
 \der u\Vert_2 \,+\, \sup_{|\msw| = h}  \: \Gamma_2(\msw,\partial u) \right ],
 \\ \displaystyle
 \Vert \partial u \Vert_2\,\left [  \Vert \partial v - 
 \der v\Vert_2 \,+\, \sup_{|\msw| = h}  \: \Gamma_2 (\msw,\partial v) \right ], 
 \end{array} \right. \end{array} \]
where $c(R)$ is $n$-independent.
\end{theorem}

The procedure by which the norm (\ref{exn2.13}) is accompanied with the norms 
(\ref{exn2.12}) is thoroughly utilized throughout this work. Therefore we define
the procedure generally. Let a collection of functions $u(R,n)$ defined by
(\ref{exn3.5}) or grid-functions ${\bf u}_n(R)$ on $G_n(R)$ be denoted by $\mathfrak{C}
(R)$. For a functional $F$ on $\mathfrak{C}(R)$ or $\mathfrak{C}(R) \times \mathfrak{C}
(R)$ the \emph{averaged value} $F_{avg}$ is defined as the arithmetic mean of functionals
over all $\mathfrak{C}(R)$ or $\mathfrak{C}(R) \times \mathfrak{C}(R)$, respectively. For
instance, the $W_2^1$-norm on $E_n(R,{\bbR}^d)$ has its averaged version denoted
by $\Vert \cdot \Vert_{avg,2,1}$ and explicitly expressed by (\ref{ex3.4}).
Another important example is a bilinear form $a(\cdot,\cdot)$ on 
$E_n(R,{\bbR}^d) \times E_n(R,{\bbR}^d)$. Its averaged version is defined by:
\[ a_{avg}(\hat{v}(n),\hat{u}(n)) \ = \ \frac{1}{vol(R)}\:\sum_R\:
 a(\hat{v}(R,n),\hat{u}(R,n)).\]

\section{Discretizations of elliptic operators}\label{sec4}
Most of constructions and descriptions in this section are given for the 
differential operator $A_0$. At the end of section we briefly explain
how to include discretizations for the lower order terms in (\ref{exp2.2}).

\subsection{Constant diffusion tensor}
Schemes and some of results regarding the convergence can be easily described
for the case of classical second order elliptic operator $A_0 = - \sum_{ij}
a_{ij} \partial_i \partial_j$, where the functions $a_{ij}$ are H\"{o}lder 
continuous on ${\bbR}^d$. At the beginning we analyze the case of constant 
diffusion tensor $\{a_{ij}\}_{11}^{dd}$.

A standard approach to a generation of discretizations of differential
operators is based on utilization of finite difference operators approximating
$\partial_i, \partial_i\partial_j$. So, by using the forward and backward difference 
operators $\der_i(p_ih_i), \widehat{\der}_i(p_ih_i)$ we can define various 
discretizations $\der_{ii}  f (\msx)$ of $\partial_i \partial_j f$.
In this way, the quadratic operator $- a_{ii} (\partial_i)^2 - a_{jj}(\partial_j)^2$ 
is approximated by the standard central difference scheme $a_{ii} \der_{ii} + a_{jj}
\der_{jj}$  with positive diagonal entries and non-positive off-diagonal entries.
For $i \ne j$ finite differences are defined by:
\begin{equation}\begin{array}{c}
\partial_i \partial_j f(\msx) \rightarrow  \der_{ij} f (\msx)= \\
 \displaystyle \frac{1}{h_1 h_2} \left \{  
 \begin{array}{c}
 f(\msx \pm \mse_i h_i \pm \mse_j h_j) - f(\msx \pm \mse_i h_i) 
 -f(\msx \pm \mse_j h_j) + f(\msx), \\
 -f(\msx \pm \mse_i h_i \mp \mse_j h_j) + f(\msx \pm \mse_i h_i) +
 f(\msx \mp \mse_j h_j) - f(\msx). \end{array}   \right.\end{array}
\end{equation}
If $a_{ij} \geq 0$, then $a_{ij} \partial_i \partial_j$
is approximated by the half sum of the first
two possibilities, otherwise by the half sum of the second two possibilities.
Hence, for the $d$-dimensional case
\begin{equation}\label{exp4.5}\begin{array}{lll} \displaystyle
 A_{\mskd \mskd \pm r_i\msed_i} &=& -\frac{1}{h^2 r_i}\:
 \Big[\frac{1}{r_i}\,a_{ii} \:-\:\sum_{m \ne i}\: \frac{1}{r_m}\,|a_{mi}|\Big],\\
  \displaystyle
 A_{\mskd \mskd \pm  (r_i\msed_i+r_j\msed_j)} &=&
 -\frac{1}{h^2 r_i r_j}|a_{ij}|, \qquad a_{ij} \geq 0,\\
  \displaystyle
 A_{\mskd \mskd \pm  (r_i\msed_i-r_j\msed_j)} &=&
 -\frac{1}{h^2 r_i r_j}|a_{ij}|, \qquad a_{ij} \leq 0,  \end{array}
\end{equation}
while the diagonal entries are calculated as the negative sum of off-diagonal ones.
To a given diffusion tensor $a = \{a_{ij}\}_{11}^{dd}$ we associate an auxiliary
tensor $\hat{a}$ defined by:
\[ \hat{a}_{ii} \ = \ a_{ii}, \quad \hat{a}_{ij} \ = \ -|a_{ij}| \quad i \ne j.\]
By using Perron-Frobenius theorem one can easily prove the following assertion.
There exist grid-steps $h$ such that the matrix $A_n$ defined 
by (\ref{exp4.5}) has the compartmental structure iff the auxiliary diffusion 
tensor $\hat{a}$ is positive definite.

In case of $d = 2$ the tensors $a, \hat{a}$ are simultaneously positive definite
or not. For higher dimensions $\hat{a}$ can be indefinite although $a$ is positive
definite. 

Apart from this standard approach to discretizations of differential operators
there can be used constructions avoiding finite difference operators. One of 
possibilities is based on the following simple result:

\begin{definition}\label{def4.1} A matrix $A_n$ with a finite band is called a 
discretization of the differential operator $A_0 = - \sum_{ij}a_{ij} \partial_i 
\partial_j$ if the equalities
\begin{equation}\label{exp4.4}
 (Au)(\msk) \ = \ \sum_{\msld} (A_n)_{\mskd \msld}\, u_{ \msld}  
\end{equation}
are valid for any polynomial $\msx \to u(\msx)$ of the second degree.
\end{definition}

\begin{lemma}\label{lem4.1}
Let $A_n$ be discretizations of the differential operator $A_0$ and
${\bf u}_n = R(\lambda,A_n) \msmu_n$ such that $\lnorm R(\lambda,A_n)\lnorm_\infty$
is bounded uniformly with respect to $n \in {\bbN}$. Then the sequence of functions
$\mathfrak{U} = \{ u(n): n \in {\bbN}\} \subset \cup_n E_n(D)$ converges in 
$\dot{C}^{(\alpha)}({\bbR}^d)$ to the solution $u = R(\lambda,A_0)f$ of (\ref{exp2.3}).
\end{lemma}

{\Proof} One can utilize arguments in standard proofs of the factor convergence 
of grid solutions as described in \cite{RM}. In order to prove the
convergence in $C^{(\alpha)}(\overline{D})$ it is necessary to have the
condition (\ref{exp4.4}) as demonstrated in \cite{LR1}. {\QED}

A useful application of this lemma is the discretization defined by:
\begin{equation}\label{exp4.6}\begin{array}{lll} \displaystyle
 A_{\mskd \mskd \pm \msed_i} &=& -\frac{1}{h^2}\:
 \Big[a_{ii} \:-\:\sum_{m \ne i}\: \frac{r_i}{r_m}\,|a_{mi}|\Big],\\
  \displaystyle 
 A_{\mskd \mskd \pm (r_i\msed_i+r_j\msed_j)} &=& -\frac{1}{h^2r_i r_j}|a_{ij}|,
 \qquad a_{ij} \geq 0, \\
  \displaystyle 
 A_{\mskd \mskd \pm (r_i\msed_i-r_j\msed_j)} &=& -\frac{1}{h^2r_i r_j}|a_{ij}|,
 \qquad a_{ij} \leq 0, \end{array}
\end{equation}
where $r_i, r_j \in {\bbN}$. Of course, a convex combination of the system matrices
(\ref{exp4.6}) is again a system matrix discretizing the classical differential operator 
$A_0 = - \sum_{ij} a_{ij} \partial_i \partial_j$. 

In the case of two dimensions the obtained structures 
of system matrix can be classified into two groups by using the 
corresponding numerical neighborhoods. The grid $G_n = \{h(k\mse_1 + l\mse_2):
k,l \in {\bbZ}\}$ has the corresponding index set of indices $\msk = (k,l)$. Possible 
numerical neighborhoods ${\cal N}(\msx)$ for the respective methods (\ref{exp4.5}), 
(\ref{exp4.6}), are illustrated in Figure~\ref{fgn4.1}.

\begin{Figure}\label{fgn4.1}
\begin{center}
\beginpicture
\setcoordinatesystem units <0.6pt,0.6pt>
\setplotarea x from 0 to 540, y from -30 to 45 
\setlinear
\plot
 10 40
 70 40 /
\plot
 90 40
150 40 /
\plot
 40 10
 40 70 /
\plot
120 10
120 70 /
\plot
 10 10
 70 70 /
\plot
 90 70
150 10 /
\plot 220 40, 280 40 /
\plot 160 10, 340 70 /
\plot 250 10, 250 70 /
\plot 360 70, 540 10 /
\plot 420 40, 480 40 /
\plot 450 10, 450 70 /
\put{$\bullet$} [c] at 10 40
\put{$\bullet$} [c] at 70 40
\put{$\bullet$} [c] at 90 40
\put{$\bullet$} [c] at 150 40
\put{$\bullet$} [c] at 40 10
\put{$\bullet$} [c] at 40 70
\put{$\bullet$} [c] at 120 10
\put{$\bullet$} [c] at 120 70
\put{$\bullet$} [c] at 10 10
\put{$\bullet$} [c] at 70 70
\put{$\bullet$} [c] at 90 70
\put{$\bullet$} [c] at 150 10
\put{$\bullet$} [c] at 40 40
\put{$\bullet$} [c] at 120 40
\put{$\bullet$} [c] at 220 40 
\put{$\bullet$} [c] at 280 40 
\put{$\bullet$} [c] at 160 10 
\put{$\bullet$} [c] at 340 70 
\put{$\bullet$} [c] at 250 10 
\put{$\bullet$} [c] at 250 70 
\put{$\bullet$} [c] at 360 70
\put{$\bullet$} [c] at 540 10 
\put{$\bullet$} [c] at 420 40
\put{$\bullet$} [c] at 480 40 
\put{$\bullet$} [c] at 450 10
\put{$\bullet$} [c] at 450 70 
\put{\emph{(a): \ $a_{12} \geq 0$, (b): \ $a_{12} \leq 0$}} [c] at 80 -10
\put{\emph{(c): \ $a_{12} \geq 0$}} [c] at 250 -10
\put{\emph{(d): \ $a_{12} \leq 0$}} [c] at 450 -10
\put{$r_1 = r_2 = 1$} [c] at 80 -30
\put{$r_1=3, r_2 = 1$} [c] at 350 -30
\endpicture
\end{center}
\end{Figure}

In order to discretize (\ref{exp2.3}) by monotone schemes we must have discretizations
$A_n$ of (\ref{exp2.2}) in terms of matrices of positive type. If we utilize 
discretizations~(\ref{exp4.6}) the following inequalities must be satisfied:
\[ \sum_{m \ne i}\: \frac{1}{hr_m}\,|a_{mi}| \ \leq \ \frac{1}{hr_i}\,a_{ii}, 
 \quad i \in \{1,2,\ldots,d\},\]
and the following Lemma follows easily.

\begin{lemma}\label{lem4.2} There exist natural numbers $r_i$ such that
the matrix $A_n$ defined by (\ref{exp4.6}) is of positive type 
iff the auxiliary diffusion tensor $\hat{a}$ is positive definite.
\end{lemma}

As already mentioned, in higher dimensions $\hat{a}$ can be indefinite although 
$a$ is positive definite. In such cases the simple schemes (\ref{exp4.6}) cannot 
give us matrices  of positive type. Rotations of coordinates are necessary to get 
system matrices of positive type. Naturally, so obtained schemes are also covered 
by Definition \ref{def4.1}.

All the results of this subsection can be straightforwardly extended to cases
with a general tensor-valued function $\msx \to a(\msx)$. Numbers $a_{ij}$ in
Expressions (\ref{exp4.5}), (\ref{exp4.6}) have to be replaced with the values
$a_{ij}(\msx)$ at the considered grid-knot $\msx \in G_n$. The condition of Lemma
\ref{lem4.2} must be satisfied uniformly with respect to $\msx \in G_n$.

\subsection{Operators in divergence form}
Discretizations of $A_0(\msx)$ which are analyzed here are constructed by using 
several rules. Rules can be described vaguely as follows. The set ${\bbR}^d$ is
discretized by $G_n(A) \subset G_n$, $u$ by ${\bf u}_n$,
$\partial_i u$ by $(U_i(p_i){\bf u}_n)$, where $p_i$ are generally
$\msx$-dependent, $a_{ij}$ are discretized by $a_{ij}(\msx_{ij}(n,\msx)),
\msx \in G_n(A)$, where $\msx_{ij}(n,\cdot): {\bbR}^d \mapsto {\bbR}^d$ are
certain functions, and integral is replaced by the 
corresponding Riemmanian sum of factors calculated at grid-knots. 
Thus the form (\ref{exp2.3}) is discretized by a sequence of forms 
${\bf v}_n, {\bf u}_n \mapsto h(n)^d \lev {\bf v}_n |A_n {\bf u}_n\des$. 
Without loss of generality, we may assume that entries
of resulting matrices $A_n$ are defined by functions $\msx_{ij}(n,\cdot)
: {\bbR}^d \mapsto {\bbR}^d$ and values $a_{ij}(\msx_{ij}(n,\msz)), \msz \in G_n$.

\subsubsection*{Standard schemes in two dimension}
We consider a subset of grid-knots $G_n(R) \ = \ \msr_0+\Big\{  k\,hr_1\,\mse_1 \:+\:l
\,hr_2\,\mse_2 \::\: k,l \in {\bbZ} \Big \} \subset G_n$ and the forms
\[ a_n(v,u) \ = \ \sum_{i\,j \,=1}^2\:\sum_{\msxd \in G_n(R)}
 \: \big(\der_i(hr_i)v\big)(\msx) \:a_{ij}\Big(\msx +\frac{h}{2}(r_1 \mse_1+
 r_2\mse_2)\Big)\: \big( \der_j(hr_j) u\big) (\msx). \]
The non-trivial off-diagonal matrix entries of $A_n$ are defined as follows:
\begin{equation}\label{exp4.1}\begin{array}{l}
 \big( A_n \big )_{\mskd \mskd\pm r_1 \msed_1} \ = \ -\: \frac{1}{h^2 r_1}
 \:\left (\,\frac{1}{r_1}\, a_{11}^{\pm +} \,-\, \frac{1}{r_2}\,|a_{12}^{\pm +}|
 \right ),\\ 
 \big( A_n \big )_{\mskd \mskd\pm r_2 \msed_2} \ = \ -\: \frac{1}{h^2 r_2}
 \:\left (\,\frac{1}{r_2}\, a_{22}^{+ \pm } \,-\, \frac{1}{r_1}\,|a_{12}^{+\pm}|
 \right ), \\
 \big( A_n \big )_{\mskd \mskd\pm (r_1 \msed_1 -  r_2 \msed_2)}
  \ = \ -\: \frac{1}{h^2 r_1 r_2}\:|a_{12}^{\pm \mp}|,\\ \end{array}
\end{equation}
and
\[\displaystyle
 a_{ij}^{\pm +} \ = \  a_{ij}\Big( \msx \:+\: \frac{h}{2} (\pm\, r_1 \mse_1 \,+\,
 r_2 \mse_2) \Big), \quad a_{ij}^{\pm -} \ = \  a_{ij}\Big( \msx \:+\:
 \frac{h}{2} (\pm\, r_1 \mse_1 \,-\, r_2 \mse_2) \Big). \]
If $a_{12} \leq 0$ and brackets in (\ref{exp4.1}) have positive values,
then the scheme (\ref{exp4.1}) defines a matrix $A_n$ with the compartmental 
structure; this scheme is called the first scheme.
Numerical neighborhoods are defined by ${\cal N}^{(-)}(\msx)
= \msx + {\cal N}^{(-)}({\bf 0})$, where
\[ {\cal N}^{(-)}({\bf 0}) \ = \ \{{\bf 0}\} \:\cup\:\Big\{ \pm
 \,hr_1\,\mse_1 \:\pm\:hr_2\,\mse_2 , \quad
 \pm\,(hr_1\,\mse_1 \:-\:hr_2\,\mse_2) \Big\}. \]
An illustration of a set ${\cal N}^{(-)}({\bf 0})$ is given in Figure 
\ref{fgn4.1}(b).

For the second scheme the grid $G_n(R)$ is defined as before while
\[\begin{array}{lll} \displaystyle
 a_n(v,u) &=& \sum_{i\,j \,=1}^2\:\sum_{\msx \in G_n(R)}
 \: \big(\der_i((-)^{i-1}hr_i)v\big)(\msx) \\ \displaystyle
 &\times & a_{ij}( \msx +\frac{1}{2}(hr_1 \mse_1 -
 hr_2\mse_2))\: \big( \der_j((-)^{j-1}hr_j) u\big)(\msx). \end{array}\]
Off-diagonal matrix entries are:
\begin{equation}\label{exp4.7}\begin{array}{l}
 \big( A_n \big )_{\mskd \mskd\pm
 r_1 \msed_1} \ = \ -\: \frac{1}{h^2 r_1} \:\left (\,\frac{1}{r_1}\,
 a_{11}^{\pm -} \,-\, \frac{1}{r_2}\,|a_{12}^{\pm -}| \right ), \\
 \big( A_n \big )_{\mskd \mskd\pm r_2 \msed_2} \ = \ -\: \frac{1}{h^2 r_2}
 \:\left (\,\frac{1}{r_2}\, a_{22}^{+\pm} \,-\, \frac{1}{r_1}\,|a_{12}^{+\pm}|
 \right ), \\
 \big( A_n \big )_{\mskd \mskd\pm (r_1 \msed_1 +  r_2 \msed_2)}
  \ = \ -\: \frac{1}{h^2 r_1 r_2}\:a_{12}^{\pm \pm}. \end{array}
\end{equation}
If $a_{12} \geq 0$ and brackets in (\ref{exp4.7}) have positive values
the compartmental structure is again attained. The set $G_n(R)$ is unchanged while
the numerical neighborhoods have different structure, ${\cal N}^{(+)}(\msx)
= \msx + {\cal N}^{(+)}({\bf 0})$, where now
\[ {\cal N}^{(+)}({\bf 0}) \ = \ \{{\bf 0}\}\:\cup\:\Big\{ \pm
 \,hr_1\,\mse_1 \:\pm\:hr_2\,\mse_2 , \quad \pm\,(hr_1\,\mse_1 \:+\:hr_2\,\mse_2) \Big\}, \]
as illustrated in Figure \ref{fgn4.1}(a). 

For each scheme there exists only one $\msx^{(\mp)}(n,\cdot)$ (instead of 3 
mappings $\msx_{ij}(n,\cdot)$) which is implicitly defined by expressions of 
$a_{ij}^{\mp \mp}$. Thus, for the first and second schemes 
we have the following respective expressions:
\begin{equation}\label{exp4.8}\begin{array}{l} \displaystyle
 \msx^{(-)}(n,\msx,\msr) \ = \ \msx \:+\:h2^{-1}\:\big( r_1\mse_1+ r_2\mse_2\big),\\
 \displaystyle
 \msx^{(+)}(n,\msx,\msr) \ = \ \msx \:+\:h2^{-1}\:\big( r_1\mse_1- r_2\mse_2\big),
 \end{array}
\end{equation}
where $\msx^{(\mp)}$ contain explicitly the parameter $\msr$.

\begin{wrapfigure}[11]{l}{200pt}
\begin{minipage}{200pt}
\begin{Figure}\label{fgn4.2}
\begin{center}
\beginpicture
\setcoordinatesystem units <1.5pt,1.5pt>
\setplotarea x from -10 to 110, y from 0 to 70 
\setlinear
\plot
-10 50
100 50 /
\plot 10 65, 35 65 /
\plot 13 67, 10 65, 13 63 /
\plot 32 67, 35 65, 32 62 /
\plot 10 45, 30 45 /
\plot 13 47, 10 45, 13 43 /
\plot 27 47, 30 45, 27 43 /
\put{$\bullet$} [c] at 10 70
\put{$\bullet$} [c] at 35 70
\put{$\bullet$} [c] at 60 70
\put{$\bullet$} [c] at 10 60
\put{$\bullet$} [c] at 35 60
\put{$\bullet$} [c] at 60 60
\put{$\bullet$} [c] at 10 50
\put{$\bullet$} [c] at 35 50
\put{$\bullet$} [c] at 60 50
\put{$\bullet$} [c] at 30 50
\put{$\bullet$} [c] at 50 50
\put{$\bullet$} [c] at 70 50
\put{$\bullet$} [c] at 10 40
\put{$\bullet$} [c] at 30 40
\put{$\bullet$} [c] at 50 40
\put{$\bullet$} [c] at 70 40
\put{$\bullet$} [c] at 10 30
\put{$\bullet$} [c] at 30 30
\put{$\bullet$} [c] at 50 30
\put{$\bullet$} [c] at 70 30
\put{$bnd(G_n(1)) \ne bnd(G_n(2)$} [c] at 40 20 
\put{$x_1 \geq 0$} [l] at 80 67 
\put{$cls(G_n(1))$} [l] at 80 60 
\put{$x_1 < 0$} [l] at 80 40 
\put{$cls(G_n(2))$} [l] at 80 33 
\put{$hr_1(1)$} [c] at 22 60 
\put{$hr_1(2)$} [c] at 20 40 
\endpicture
\end{center}
\end{Figure}
\end{minipage}
\end{wrapfigure}
\par
In a general two-dimensional case the function $a_{12}$ may vary and
the compartmental structure of $A_n$ can be broken at some grid-knots.
Besides, the function $a_{12}$ changes sign on ${\bbR}^2$ and the compartmental 
structure cannot be realized by using only one scheme for all the grid-knots.
In the case of continuous diffusion tensor one can prove the following result.
There exists a partition ${\bbR}^d = \cup_{\mp}\cup_l D_l^{(\mp)}$ such that
the function $a_{12}$ does not change sign on $D_l^{(\mp)}$. To each
$D_l^{(\mp)}$ there correspond grid-steps $hr_i^{(\mp)}(l)$
such that the sets $G_n^{(\mp)}(l) = D_l^{(\mp)} \cap G_n$ are subgrids and
$G_n(R) = \cup_{\mp l}G_n^{(\mp)}(l)$. To each grid-knot in $int(G_n^{(\mp)}(l))$
we associate one of neighborhoods ${\cal N}^{(\mp)}(\msx)$ and entries of $A_n$
defined by (\ref{exp4.1}) or (\ref{exp4.7}). Construction of entries of
$A_n$ at grid-knots on the boundary of sets $D_l^{(\mp)}$ can bear
a considerable amount of difficulties. For instance, one of simple cases is 
given by $G_n = G_n(1) \cup G_n(2)$ as illustrated in Figure \ref{fgn4.2},
where we are faced with difficulties regarding the constructions of sets 
${\cal N}(\msx), \msx \in \{{\bbR}^2: x_2=0\}$ and the corresponding entries of $A_n$.

\subsubsection*{Extended schemes in two dimensions}
The mappings $\msx^{(\mp)}(n,\msx,\msr)$ in ${\bbR}^2$ are defined by (\ref{exp4.8}), 
where the superscripts 
$(\pm)$ are related to the sign of $a_{12}$, {\em i.e.} to the cases $a_{12} \geq 0$ and
$a_{12} \leq 0$, respectively. The bilinear form on $G_n \times G_n$
\begin{equation}\label{exp4.9}\begin{array}{c} \displaystyle
 a_n^{(-)}(v,u) \ = \ \sum_{\msxd \in G_n}\: \left (\Big[\sum_{i = 1}^2\:
 a_{ii}(\msx^{(-)}(n,\msx))\, \big(\der_i(h)v \big)(\msx) \big(\der_i(h)u \big)(\msx)
 \right .\\  \displaystyle
 \:+\:\sum_{i \ne j}\:a_{ij}(\msx^{(-)}(n,\msx,\msr))\, \big(\der_i(r_ih)v \big)(\msx)
 \big(\der_j(r_jh)u \big)(\msx)\Big] \\  \displaystyle \left .
 +\: \sum_{i \ne j}\:a_{ij}(\msx^{(-)}(n,\msx,\msr))\, \frac{r_i}{r_j}\:
 \Big[ \big(\der_i(h)v \big)(\msx) \big(\der_i(h)u \big)(\msx)
 \:-\:\big(\der_i(r_ih)v \big)(\msx) \big(\der_i(r_ih)u \big)(\msx)\Big]\right )
 \end{array}
\end{equation}
can be written as
\[ a_n^{(-)}(v,u) \ = \ \lev {\bf v}\,|\,A_n\,{\bf u}\des,\]
where $\msr \in {\bbN}^2$, $\msx^{(\mp)}(n,\msx) =\msx^{(\mp)}(n,\msx,{\bf 1})$ by
convention, ${\bf v}, {\bf u}$ are columns with entries $v(\msx), u(\msx)$, and
$A_n$ is the resulting system matrix. 
By utilizing the standard variational calculus we can 
derive entries of $A_n$. In order to rewrite the entries of $A_n$ in a concise form we 
use the following abbreviations:
\[\begin{array}{l}\displaystyle
 a_{ij}^{\pm +}(\msr) \ = \  a_{ij}\big( \msx \:+\: 2^{-1} (\pm\, r_1h \mse_1 \,+\,
 r_2h \mse_2) \big),\\ \displaystyle
 a_{ij}^{\pm -}(\msr) \ = \  a_{ij}\big( \msx \:+\:
 2^{-1} (\pm\,r_1 h \mse_1 \,-\,r_2 h \mse_2) \big), \end{array}\]
and $a_{ij}^{\alpha \beta} = a_{ij}^{\alpha \beta}({\bf 1})$ by convention.
For instance, $a_{ij}^{+ +}(\msr) = a_{ij}(x^{(+)}(n,\msx,\msr)),
a_{ij}^{- +}(\msr) = a_{ij}(x^{(+)}(n,\msx,\msr)-r_1h \mse_1)$, {\em etc.}
Now we have the following result:
\begin{equation}\label{exp4.10}\begin{array}{l}
 \big( A_n \big )_{\mskd \mskd\pm  \msed_1} \ = \ -\: \frac{1}{h^2}
 \:\left (a_{11}^{\pm +} \,-\, \frac{r_1}{r_2}\,|a_{12}^{\pm +}
 (\msr)| \right ),\\ 
 \big( A_n \big )_{\mskd \mskd\pm \msed_2} \ = \ -\: \frac{1}{h^2}
 \:\left (a_{22}^{+ \pm } \,-\, \frac{r_2}{r_1}\,|a_{12}^{+\pm}
 (\msr)| \right ), \\
 \big( A_n \big )_{\mskd \mskd\pm (r_1 \msed_1 - r_2  \msed_2)}
  \ = \ -\: \frac{1}{r_1 r_2 h^2}\:a_{12}^{\pm \mp}(\msr).
 \end{array}
\end{equation}
Obviously, the constructed schemes are analogs of the schemes (\ref{exp4.6}) which
are derived for the classical elliptic differential operator $A_0 = -\sum_{ij} a_{ij}
\partial_i \partial_j$. 
If $a_{12} \leq 0$ and 
\[ \inf_{\msxd \in G_n}\: \left\{
 \frac{1}{r_ih}\, a_{ii}^{\alpha \beta } \,-\, \frac{1}{r_jh}\,|a_{12}^{\gamma \delta}
 (\msr)|\right\} \ > \ 0, \qquad \alpha, \beta, \gamma, \delta \:\in\: \{=,-\},\]
the constructed matrix $A_n$ is compartmental. We call it the first extended scheme.

The second extended scheme can be constructed analogously starting from another
bilinear form. Let us define $a_n^{(+)}(v,u)$ form $a_n^{(-)}(v,u)$ by the following
simple change. Each $\der_i(r_ih)$ in (\ref{exp4.9}) has to be replaced with $
\der_i((-)^{i-1}hr_i)$. The
obtained form is denoted by $a_n^{(+)}(v,u)$ and the corresponding system matrix
$A_n$ has entries:
\begin{equation}\label{exp4.11}\begin{array}{l}
 \big( A_n \big )_{\mskd \mskd\pm
 \mse_1} \ = \ -\: \frac{1}{h^2} \:\left (
 a_{11}^{\pm -} \,-\, \frac{r_1}{r_2}\,|a_{12}^{\pm -}(\msr)| \right ), \\
 \big( A_n \big )_{\mskd \mskd\pm \mse_2} \ = \ -\: \frac{1}{h^2}
 \:\left ( a_{22}^{+\pm} \,-\, \frac{r_2}{r_1}\,|a_{12}^{+\pm}(\msr)|
 \right ), \\
 \big( A_n \big )_{\mskd \mskd\pm (r_1 \msed_1 +  r_2 \msed_2)}
  \ = \ -\: \frac{1}{r_1 r_2 h^2}\:a_{12}^{\pm \pm}(\msr), \end{array}
\end{equation}
in a complete analogy with (\ref{exp4.6}). It is compartmental for $a_{12} \geq 0$.

It is important to point out that the schemes (\ref{exp4.10}), (\ref{exp4.11}) are
obtained by a guess from the schemes (\ref{exp4.6}), {\em i.e.} they are a priori defined 
from some general principle which is not based either on difference equations or
variational equalities. Then bilinear forms (\ref{exp4.9}) are constructed so that 
variational equalities are equivalent with the schemes. This non-conventional 
approach in the construction of schemes has an apparent drawback. The forms 
(\ref{exp4.9}) do not seam to be positive definite. Actually, their positive 
definiteness can be ensured in the $\lnorm \cdot\lnorm_{avg,2,1}$-norm, {\em i.e.}
a weaker norm than $\lnorm \cdot\lnorm_{2,1}$-norm. This problem is analyzed in the
next subsection.

In the case of a general function $a_{12}$ on ${\bbR}^2$ we cannot use only one of
proposed two extended scheme. Here we give a construction of $A_n$, which 
discretize the differential operator $A_0$, possessing the compartmental structure.
Our intention is to construct a sequence $a^{(r)} = \{ a_{ij}^{(r)}
\}_{11}^{22}$ of tensor-valued functions converging $\ast$-weakly in $L_\infty({\bbR}^2)$
to the original tensor-valued function~(\ref{exp2.2}), such that:
\begin{enumerate}
 \item Each $a^{(r)}$ defines a form $v,u \mapsto a^{(r)}(v,u)$ and for any fixed 
pair $v,u \in W_2^1({\bbR}^2)$
\[  a(v,u) \ = \ \lim_r \:a^{(r)}(v,u). \]
 \item The forms $a^{(r)}(\cdot,\cdot)$ can be discretized by using described extended
schemes. The resulting discretizations are denoted by $h^2a_n^{(r)}(\cdot,\cdot)$.
 \item The matrices $A_n^{(r)}$ defined by $a_n^{(r)}(v,u) = \lev {\bf v}|A_n^{(r)}
{\bf u} \des$ have the compartmental structure.
\end{enumerate}

In the first step of construction we assume that the functions $a_{ij}$ are
uniformly continuous on ${\bbR}^2$. For each $r \in {\bbN}$ we define the sets
\[\begin{array}{lll}
 F_r(+) \ = \ \{\msx \in {\bbR}^2 \::\: a_{12}(\msx)
 \: \geq \ r^{-1}\},\\
 F_r(-) \ = \ \{\msx \in {\bbR}^2 \::\: a_{12}(\msx)
 \: \leq \ -\,r^{-1}\},\\
 F_r(0) \ = \ {\bbR}^2 \ \setminus \ F_r(+)\,\cap\,F_r(-), \end{array}
 \qquad r \in {\bbN}.\]
Due to the uniform continuity of functions $a_{ij}$ there exists a positive 
$\varepsilon_r$ such that for sufficiently large $r$ and $n$:
\begin{enumerate}
 \item $dist(F_r(+),F_r(-)) = \varepsilon_r > 0$.
 \item The set $F_r(0)$ has a finite number of connected components.
 \item To each connected component of $F_r(0)$ there corresponds one connected
set of $F_r(0)  \cap  int \big(\cup_{\msx \in G_n} C(h,\msx) \big)$.
\end{enumerate}

Now we define the sets
\[\begin{array}{l} \displaystyle
 G_n^{(\pm)} \ = \ \cup_{\msx \in F_r(\pm)\cap G_n}\: \overline{C(h, \msx)},\\ \displaystyle
 G_n^{(0)} \ = \ {\bbR}^2 \ \setminus \ \Big( G_n^{(+)} \,\cup\,G_n^{(-)}\Big),\end{array}\]
and point out that for sufficiently large $n$ the sets $F_r(0)$ and $G_n^{(0)}$ 
have the same number of connected components. In addition, each connected component
is equal to the union of interiors of certain cubes $C(h,\msx)$ 
with vertices $\msx \in F_r(0)$. In order to complete a partition of 
${\bbR}^2$ into subsets $G_n^{(\pm)}$ we can extend the constructed sets 
$G_n^{(\pm)}$ in the following way. Each non-trivial $C(h,\msx) \cap G_n^{(0)}$ can 
be attached to the nearest existing $G_n^{(\pm)}$. The first and second extended
schemes have to be used for grid-knots in so completed $G_n^{(\pm)}$, respectively.
The associated forms $a_n^{(\pm)}(v,u)$ are defined as in (\ref{exp4.9}) with the following
difference. Instead of the sum over grid-knots $\msx \in G_n$, in the case of 
$a_n^{(\pm)}(v,u)$, we must have the sum over grid-knots in $G_n^{(\pm)}$, respectively.
Then $a_n^{(\pm)}(v,u) = \lev {\bf v}_n|A_n {\bf u}_n\des$ and $A_n$ is compartmental.
Of course, the parameter $\msr$ varies on $G_n$ so appropriate notation is
$\msr(\msx)$.

In the last step of the construction we assume that the functions $a_{ij}$ are measurable.
We have to smooth the functions $a_{ij}$ by a mollifier $\{\chi_n:n \in {\bbN}\}$,
where $\chi_n(\msx) = n^2\chi_1(n\msx)$
and $\chi_1$ is a non-negative function with support in the disc $B_1({\bf 0})$. 
In this way we get a sequence of tensor-valued functions $a(n) = a \ast \chi_n$ 
converging $\ast$-weakly in $L_\infty({\bbR}^2)$ to the original tensor-valued
function $a$. The previous steps of constructions must be applied to the
differential operators $A_0^{(r)}$ with coefficients $a_{ij}(n)$.

\subsection{Higher dimensions}
The index set of pairs $I(d) = \{ \{ij\} : i < j,\:
i,j = 1,2,\ldots,d, i \ne j\}$ has the cardinal number $m(d)=d(d-1)/2$.
To each index $\{kl\} \in I(d)$ we associate three coefficients,
\begin{equation}\label{exp4.12}
 a_{kk}^{\{kl\}} = \frac{1}{d-1}\, a_{kk}, \quad a_{ll}^{\{kl\}} =
 \frac{1}{d-1}\, a_{ll}, \quad a_{kl}^{\{kl\}} = a_{kl},
\end{equation}
and a bilinear form $a^{\{kl\}}(\cdot, \cdot)$,
\begin{equation}\label{exp4.13}
 a^{\{kl\}}(v,u) \ = \sum_{i,j \in \{r,s\}}\: \int_D \: a_{ij}^{\{kl\}}(\msx) \:
 \partial_i v(\msx)\, \partial_j u(\msx) \,d\msx.
\end{equation}
Apparently, for each $v \in \dot{W}_\infty^1({\bbR}^d), \, u \in
\dot{W}_1^1({\bbR}^d)$ with
compact supports, the following equality is valid:
\[ a(v,u) \ = \ \sum_{\{kl\} \in I(d)} \, a^{\{kl\}}(v,u) . \]
To each of the forms $a^{\{kl\}}(\cdot, \cdot)$ we can associate a sequence of
forms $a_n^{\{kl\}}(\cdot,\cdot)$ and matrices $A_n^{\{kl\}}$ constructed by using
extended schemes. Then the matrix 
\begin{equation}\label{exp4.14}
 A_n \ = \ \sum_{\{kl\} \in I} \:A_n^{\{kl\}} ,
\end{equation}
is a discretization of $A_0(\msx)$. In the construction of $A_n^{\{kl\}}$, 
pairs $\msr^{\{kl\}}(\msx) = (r_1^{\{kl\}}(\msx),r_2^{\{kl\}})(\msx) \in {\bbN}^2$ can be used,
generally depending on $\msx$. 

In the remaining part of this subsection the objects of analysis are the compartmental 
structure of (\ref{exp4.14}) and strict ellipticity of forms ${\bf v}, {\bf u} \mapsto
\lev {\bf v}|A_n {\bf u}\des$. 

For $\msx = h\msk \in G_n$ the neighborhoods ${\cal N}(\msx)$ have
structures which are mutually similar one to the other. If all $a_{ij}(\msx) \ne 0$ 
then the minimal
number of elements in ${\cal N}(\msx)$ is $1+d+d^2$. In this case the set 
${\cal N}(\msx)$ consists of its center, $2d$-grid-knots on the $d$-dimensional
cross $\{\pm \mse_i: i = 1,2,\ldots,d\}$ and 2 grid-knots in each two-dimension
plane. Generally, the number of grid-knots in two-dimensional plane may be
larger than 2. Here we consider only the case of at most two grid-knots in 
two-dimensional planes.
This demand has the following implication on the construction of discretizations
$A_n^{(rs)}$. Let the pairs $\mse_r, \mse_s$ and $\mse_s, \mse_t$ define
two-dimensional planes and let $A_n^{(rs)}, A_n^{(st)}$ be the corresponding
discretizations which are constructed by using parameters 
$\msr^{(rs)}, \msr^{(st)}$. Then there must hold $(\msr^{(rs)})_s =
(\msr^{(st)})_s$. In such case the off-diagonal entries
of $A_n$ have the following general structure:
\begin{equation}\label{exp4.15}\begin{array}{lll}
 \big(A_n\big)_{\mskd \mskd\pm \msed_i} &=&    -\frac{1}{h^2} \Big[ \overline{a}_{ii}
 \,-\,\sum_{m \ne i}\: \frac{r_i}{r_m}|a_{im}(\msx_{im})|\Big],\\
 \big(A_n\big)_{\mskd \mskd\pm \msmd_{ij}} &=& -\frac{1}{r_i r_j h^2}
 |a_{ij}(\msx_{ij})|, \end{array}
\end{equation}
where $\overline{a}_{ii}$ are certain convex combinations of values $a_{ii}(
\msx_{ii}^{\{kl\}})$, $\msx_{ij}^{\{kl\}}$ are certain points close to $\msx$ which
are defined by the extended schemes and $\msm_{ij}$ is either $r_i\mse_i + r_j\mse_j$
or $r_i\mse_i - r_j\mse_j$, depending on the sign of $a_{ij}$.

Let us assume that the auxiliary tensor-valued function $\msx \mapsto \hat{a}(
\msx) = \{\hat{a}_{ij}(\msx)\}_{11}^{dd}$ fulfils the double inequality:
\begin{equation}\label{exp4.16}
 \Md(\hat{a})\,|\msz|^2 \ \leq \ \sum_{i,j} \: z_i\,\hat{a}_{ij}(\msx)\,z_j
 \ \leq \  \Mg(\hat{a})\,|\msz|^2, \quad \msx \in {\bbR}^d.
\end{equation}

The following result is valid:

\begin{proposition}\label{Prop4.1} Let $a_{ij}$ be uniformly continuous on 
${\bbR}^d$ and let the double inequality (\ref{exp4.16}) be valid. There exist 
a function $\msx \mapsto \msr(\msx) \in {\bbN}^d$, to be used in the construction
of matrices $A_n^{\{kl\}}$, such that (\ref{exp4.14}) has the compartmental structure.
\end{proposition}

In order to prove this proposition one has to apply Perron-Frobenius theorem 
to the matrix $\hat{a}(\msx)^{-1}$ and construct the eigenvector $\msz(\msx)$,
corresponding to the minimal eigenvalue of $\hat{a}(\msx)$, with the property
$z_i(\msx) \geq \underline{m} > 0$, where $\underline{m}$ is $\msx$-independent.

If the coefficients $a_{ij}$ are not continuous, we have to use a smoothing procedure
and apply the derived results to smeared out coefficients $a_{ij}(n) =
a_{ij}\ast \chi_n$.

Let us turn to the problem of strict ellipticity.
A discrete form $a_n(\cdot,\cdot)$ on $l_0(G_n(R)) \times l_0(G_n(R))$ is said to be 
strictly elliptic~\cite{Yo} if there exist two positive numbers $\Md(a_n), \Mg(a_n)$ such
that
\[ \Md(a_n)\,\sum_{i=1}^d\:\lnorm U_i(r_i) {\bf u} \lnorm_{R2}^2 \ \leq \ 
 \lev {\bf u} \,|\,A_n\,{\bf u} \des_R \ \leq \ 
 \Mg(a_n)\,\sum_{i=1}^d\:\lnorm U_i(r_i) {\bf u} \lnorm_{R2}^2.\]
We shall say that a sequence of discrete forms $a_n(\cdot,\cdot)$ is strictly
elliptic uniformly with respect to $n \in {\bbN}$ if the positive numbers 
$\Md(a_n), \Mg(a_n)$ are independent of $n$.

The quantities $\overline{a}_{ii}$ in~(\ref{exp4.15}) have the general structure
\[\begin{array}{c}
 \omega_n(a_{ii},\msx) \ = \ \sum_{s=1}^d \: f_s(\msx)\,a_{ii}(x_{ii}(n,\msx,s)),\\
 f_s(\msx) \ \geq 0, \quad \sum_{s=1}^d \: f_s(\msx) \ = \ 1,
 \end{array}\]
where $\msx_{ij}(n,\msx,s)$ are defined by the construction of schemes.

\begin{proposition}\label{Prop4.2} Let $A_0 = -\sum \partial_i a_{ij}\partial_j$ be
discretized by matrices $A_n$ with the following properties:
\begin{description}\itemsep 0.cm
 \item{a)} Each $A_n$ has the compartmental structure whenever
the double inequality (\ref{exp4.16}) holds. 
 \item{b)} There exist $p_i, t_i \in {\bbR}$ such that
\[ \big(A_n\big)_{\mskd \mskd\pm r_i\msed_i} \ = \ -\frac{1}{h^2}\,p_i\:
 \Big[ \omega_n(a_{ii},\msx) \,-\,\sum_{m \ne i} \: t_i\,|a_{im}(\msx_{im}(n,\msx))|
 \Big] \]
at each grid-knot $\msx = h \msk$. 
 \item{c)} If $\msm$ is not in the direction of $\pm \mse_i, i = 1,2,\ldots,d$, then
$\big(A_n\big)_{\mskd \mskd+\msmd}$ does not depend on any of $a_{ii}, i = 1,2,\ldots,d$.
\end{description}

Then the discrete forms ${\bf v}, {\bf u} \mapsto
\lev {\bf v}|A_n{\bf u}\des_R$ are strictly elliptic on $l_0(G_n(R)) \times l_0(G_n
(R))$ uniformly with respect to $n \in {\bbN}$.
\end{proposition}

{\Proof} Let us consider diffusion tensors $a, b$ where $a = \{a_{ij}\}_{11}^{dd}$
and $b$ is defined by $b_{ij}=a_{ij}-\kappa \delta_{ij}$. The corresponding
auxiliary tensors are denoted by $\hat{a}, \hat{b}$ as usual. Due to the fact that
$a, \hat{a}$ are positive definite uniformly with respect to the points of ${\bbR}^d$
there can be chosen $\kappa > 0$ so that $b, \hat{b}$ are also positive
definite on ${\bbR}^d$. Let us define matrices $H_n$ on $G_n(R)$ by the following
non-trivial entries:
\[ \big(H_n\big)_{\mskd \mskd} \ = \ \frac{2\kappa }{h^2}\sum_{i=1}^d\:\frac{1}{r_i^2}, \quad
 \big( H_n \big)_{\mskd \mskd \pm r_i\msed_i} \ = \ -\frac{\kappa}{r_i^2 h^2}, \  
 i = 1,2,\ldots,d.\]
We have $\lev {\bf v}|H_n{\bf u}\des_R = \kappa \sum_{i=1}^d \lnorm U_i(r_i){\bf u}
\lnorm_{R2}^2$. The matrix $B_n = A_n - H_n$ has the compartmental structure in 
accordance with assumptions of this proposition. A matrix $B_n$ with the compartmental
structure is necessarily positive definite, {\em i.e.} there must hold
$\lev {\bf u} |B_n {\bf u}\des_R \geq 0$ for any ${\bf u} \in l_0(G_n(R))$. Therefore
\[\lev {\bf u} \,|\,A_n \,{\bf u}\des_R \ = \ \lev {\bf u} \,|\,H_n \,{\bf u}\des_R \:+\:
 \lev {\bf u} \,|\,B_n \,{\bf u}\des_R \ \geq \ \lev {\bf u} \,|\,H_n \,{\bf u}\des_R \:=\:
 \kappa \sum_{i=1}^d \lnorm U_i(r_i){\bf u}\lnorm_{R2}^2, \]
proving the assertion. {\QED}

The first order differential operators of (\ref{exp2.2}) can be discretized by using
upwinding method. In the case of continuous functions $b_j$ the resulting 
system matrices have the compartmental structure. If the coefficients $b_j$ are
not continuous, the compartmental structure is obtained for sufficiently large $n$.

\section{Convergence of numerical solutions}\label{sec5}
The difference schemes which are constructed in Section \ref{sec4} enable us to 
generate grid-solutions ${\bf u}_n \in L(G_n(R))$ and the corresponding imbedding
$u(n) =\Phi_n(R){\bf u}_n$ into $E_n(R,{\bbR}^d)$ spaces. Basic schemes are
derived from variational equalities
\begin{equation}\label{exp5.1}
 \lev {\bf v}_n\,|\,\lambda I+A_n\,{\bf u}_n \des \ = \  \lev  {\bf v}_n\,|\,
 \msmu_n\des, \quad {\bf v}_n \in w_2^1(G_n),
\end{equation}
while extended schemes follow from a general
principle which is not a priori related to any bilinear form and variational
equalities. For a differential operators in divergence form $A_0$ the derived 
discretizations $A_n$ have the compartmental structure, {\em i.e.} they satisfy 
(\ref{exp4.4}) only for constant
functions. Consequently, as expected, we cannot use the result on convergence
of Section \ref{sec4}. In order to prove the convergence of approximate solutions 
we must replace conditions of Lemma \ref{lem4.1} with some other property. It is
natural to look for a formulation in terms of variational equalities (\ref{exp5.1})
and follow well known methods of finite elements. Unfortunately, we are faced with 
a difficulty. For extended schemes the bilinear form (\ref{exp4.9}) does not seam
to be strictly elliptic. The strict ellipticity of forms is a basic supposition
in a proof of $W_2^1$-convergence. For this purpose we can use a general result of
Proposition \ref{Prop4.2}. 
The underlining idea is clear from analyzing the 2-dimensional problem with 
$a_{12} \leq 0$ on ${\bbR}^2$ and a fixed $\msr = (r_1,r_2)$. 

When Proposition \ref{Prop4.2} is applied to the present case we get

\begin{corollary}\label{cor5.1} Let the auxiliary diffusion tensor $\hat{a}$ be
strictly positive on ${\bbR}^d$. There exists a mapping $\msx \mapsto \msr(\msx)
\in {\bbN}^d$ such that the matrices $A_n$ of (\ref{exp4.14}) have the compartmental
structure, and the forms defined by $a_n(v,u) = \lev {\bf v}|A_n {\bf u}\des$ are
strictly elliptic uniformly with respect to $n \in {\bbN}$.
\end{corollary}

For extended schemes~Corollary~\ref{cor5.1} can be
expressed simply as:
\[ \Md \,\sum_{i=1}^2\,\lnorm U_i{\bf u}_n 
 \lnorm_2^2 \ \leq \ a_n(u,u) \ \leq \ \Mg \,\sum_{i=1}^2\,\lnorm U_i
 {\bf u}_n \lnorm_2^2\]
with certain $\Md, \Mg$. For basic schemes we must use averaged norms. In order
to demonstrate that the averaged norms have no effect on the convergence of
approximate solutions we proceed further by using averaged norms.

\subsection{Consistency}
For constant functions $a_{ij}$ we have an obvious equality:
\[ a(v,u) \ = \ \lim_n \:h^d\,a_{avg,n}(\hat{v}(n),\hat{u}(n))\]
implied by (ii) of Theorem \ref{thn3.1} and Corollary \ref{corn3.1}. For non-constat 
functions $a_{ij}$ this property can be proved. Actually, we need a more
general result.

\begin{proposition}[Consistency]\label{Prop5.1} Let $\mathfrak{V} = \{v(n)
: n \in {\bbN}\}$ and $\mathfrak{U} = \{u(n): n \in {\bbN}\}$,
$v(n), u(n) \in E_n({\bbR}^2)$, converge weakly and strongly in 
$W_2^1({\bbR}^d)$ to $v, u \in W_2^1({\bbR}^d)$, respectively. Then
\[ a(v,u) \ = \ \lim_n \: h^d\,\lev {\bf v}_n\,|\,A_n \,{\bf u}_n\des.\]
\end{proposition}

The present object of analysis is the functional $F_n(\cdot,\cdot)$ on 
$E_n(R,{\bbR}^d) \times E_n(R,{\bbR}^d)$, defined by
\[ F_n(v,u) \ = \ h^d\,vol(R)\,\sum_{\mskd \in I_n(R)} \:p(\msx(n,h\msk))\,
 \der_i(r_ih) v(h\msk)\, \der_j(r_jh) u(h\msk),\]
where $p$ is a bounded and uniformly continuous function on ${\bbR}^d$,
$\msx(n,\cdot)$ is a transformation of ${\bbR}^d$ such that
$| \msx(n,h\msk))-h\msk| \leq \omega n^{-1}$, $\omega$ is fixed.
Obviously, the form $h^d a_n(\cdot,\cdot)$ of (\ref{exp4.9}) is a finite sum of 
$F_{avg,n}$ with various $p = a_{ij}$. Therefore a proof of Proposition follows from
the following result.

\begin{lemma}\label{lem5.4} Let the sequence 
$\mathfrak{V} = \{ v(n): n \in {\bbN}\}\subset \cup_n E_n(R,{\bbR}^d)$ 
converge $W_2^1$-weakly to $v \in W_2^1({\bbR}^d)$ and the sequence 
$\mathfrak{U} = \{ u(n): n \in {\bbN}\}\subset \cup_n E_n(R,{\bbR}^d)$ converge 
$W_2^1$-strongly to $u \in W_2^1({\bbR}^d)$. Then
\[ \lim_n \:F_n(v(n),u(n)) \ = \ \int\: p(\msx)\,\frac{\partial v}{\partial x_i}(
 \msx)\,\frac{\partial u}{\partial x_j}(\msx)\:d\msx.\]
\end{lemma}

{\Proof} The integral on the right hand side is denoted by $P(v,u)$. We have
to analyze $P(v(n),u(n))$ as $n \to \infty$. To simplify notation we assume $i, j \in
\{1,2\}$. The indices $\msk = (k_1,k_2,\ldots,k_s)$ are denoted shortly as $\msk = (k_1,
\msk')$ and $\msk = (k_1,k_2,\msk'')$, $k_1,k_2 \in {\bbZ}$. The corresponding
$r_1, r_2$ are denoted by $r,s$, respectively. After inserting expressions for
$v(n), u(n)$ and carrying out a straightforward calculation,
we get the following expression for the case of $i = j =1$:
\[ P(v(n),u(n)) \ = \ \sum_{k,\mskd',\msld'\in I_n(R)}\: \int_{J(k,r)}\:dz\:
 (\psi_{\mskd'}|p(z,\cdot) \psi_{\msld'}) \: v_{(k\mskd') i}\,u_{(k \msld') i},\]
where $v_{\mskd i}=\der_i(r_ih) v(h\msk), u_{\mskd j}=\der_j(r_jh) u(h\msk)$ and
$J(k,r) = [kh,kh+rh]$. Analogously we get for $i=1, j=2$:
\[\begin{array}{c} \displaystyle
 P(v(n),u(n)) \ = \ \sum_{k,l,\mskd'',\msld''\in I_n(R)}\: \int_{J(k,r) \times
 J(l,s)}\:dz_1 dz_2\: (\psi_{\mskd''}|p(z_1,z_2,\cdot) \psi_{\msld''}) \\
 \displaystyle \times \:
 \Big[v_{(k,l,\mskd'') i}\,u_{(k,l, \msld'') j}\,\psi_k(z_1)\psi_l(z_2) \,+\,
 v_{(k,l,\mskd'') i}\,u_{(k+r,l, \msld'') j}\,\psi_{k+r}(z_1)\psi_l(z_2) \\
 \displaystyle  +\:
 v_{(k,l+s,\mskd'') i}\,u_{(k,l, \msld'') j}\,\psi_k(z_1)\psi_{l+s}(z_2) \,+\,
 v_{(k,l+s,\mskd'') i}\,u_{(k+r,l, \msld'') j}\,\psi_{k+r}(z_1)\psi_{l+s}(z_2)
 \Big]. \end{array}\]
Now we approximate the obtained expressions by inserting specific values of $p$
in the integrals over intervals $J(k,r)$. The error arised from approximation can be 
estimated by a bound depending on the product of 
$osc(p,n) = \sup \{|p(\msx)-p(\msy)|: \msx, \msy \in {\bbR}^d, |x_i - y_i|
\leq 4r_ih\}$ and the norms $\lnorm U_i(r_i){\bf v}\lnorm_{R2}, \lnorm U_j(r_j){\bf u}
\lnorm_{R2}$. For instance, for $i=j=1$ the resulting expressing is
\[ Q_n(v(n),u(n)) \ = \ h^d\,vol(R)\,\sum_{k,\mskd', \msld' \in I_n(R)} \:
 p(\msx(n,h\msk))\, s_{\mskd' \msld'} \,v_{(k,\mskd') i}\,u_{(k,\msld') j},\]
where generally $s_{\mskd \msld} = (\psi_{\mskd}|\psi_{\msld})\Vert \psi_{\mskd}
\Vert_1^{-1}$. The quantities $s_{\mskd \msld}$ have properties $s_{\mskd \msld} \geq 0, 
\sum_{\msld} s_{\mskd \msld} = 1$. In addition
\[\begin{array}{c}
 |P(v(n),u(n)) \,-\,Q_n(v(n),u(n))| \ \leq \ osc(p,n)\,h^d\,
 \lnorm {\bf v}\lnorm_{R2}\, \lnorm {\bf u}\lnorm_{R2} \\
 \leq \ osc(p,n)\, (1-\sigma^2)^{-1}\,\Vert v(n)\Vert_{2,1}\, \Vert u(n)\Vert_{2,1}.
 \end{array}\]
Instead of $P(v(n),u(n))$ we consider in the following $Q_n(v(n),u(n))$.
The quantity $Q_n$ would be equal $F_n$ if $s_{\mskd \msld}$ were absent and the double
sum were replaced with the single sum over indices $\msk$. 
\[ Q_n(v(n),u(n))\,-\,F_n(v(n),u(n)) \ = \
 h^d\,vol(R)\,\sum_{\mskd \msld} \:p(\msx(n,h\msk))\,v_{\mskd i}\,\delta_{kl}\,
 s_{\mskd' \msld'} \,\big(u_{\msld j}-u_{\mskd j}\big).\]
To estimate the right hand side we need $\overline{p} = \sup p$.
\[\begin{array}{c} \displaystyle
 |Q_n(v(n),u(n))\,-\,F_n(v(n),u(n))| \ \leq \\ \displaystyle
 \overline{p}\:h^{d/2}\lnorm U_i(r_i) {\bf v}_n \lnorm_{R2} \:h^{d/2}\,
 \max \{\lnorm \big(Z(r_j,j)-I\big) U_j(r_j){\bf u}_{n}\lnorm_{R2} \,:\,
 j = 1,2,\ldots,d\}. \end{array}\]
Due to the strong convergence of $\mathfrak{U}$, Lemma \ref{lem3.2} and Corollary
\ref{corn3.1} we have for $\msw = hr_j \mse_j$:
\[ h^{d/2}\,\lnorm \big(Z(r_j,j)-I\big) U_j(r_j){\bf u}_{n}\lnorm_{R2} 
 \ \leq \ (1-\sigma^2)^{-1/2}\:\Vert \big(Z(\msw)-I\big) \,u(n)\Vert_2 \]
so that $\lim_n F_n(v(n),u(n)) = \lim_n P(v(n),u(n)) = P(v,u)$. {\QED}

\subsection{$W_2^1$-convergence}
It is described until now how the form $\lambda (v|u)+a(v,u)$ is discretized by forms 
$h^d \lambda  \lev {\bf v}_n|{\bf u}_n\des+\lev {\bf v}_n|A_n{\bf u}_n\des$. 
In order to solve the discretized problem (\ref{ex2.10}) or (\ref{exp5.1}), we
have to describe how to discretize the linear function $v \to \lev v|\mu\des$
by $h^d\lev {\bf v}_n|\msmu_n\des$, where $\msmu_n \in l(G_n(R))$.
First we have to demonstrate the existence of $\msmu_n$ such that
$h^d\lev {\bf v}_n|\msmu_n\des \to \lev v|\mu\des$.

\begin{lemma}\label{lem5.2} Let $\mu$ be a continuous linear functional on 
$W_2^1({\bbR}^d)$. There exists discretizations $\msmu_n(R) \in l(G_n(R))$ such that
\begin{equation}\label{exp5.3}
 \lim_n\: \left |\frac{ h^d\,\lev {\bf u}_n(R) \,| \,\msmu_n(R)\des_R \,-\, 
 \lev u(R,n)\,|\,\mu \des}{\Vert u(R,n) \Vert_{2,1}}\right | \ = \ 0
\end{equation}
for any $R$.
\end{lemma}

{\Proof} It suffices to consider the case $\mu = \partial f, f \in L_2({\bbR}^d)$,
$\Vert \mu \Vert_{2,-1} = \Vert f \Vert_2$. Let us define
\[ \big(\msmu_n \big)_{\mskd} \ = \  -\:h^{-d}\,(\partial \psi_{\mskd}|f).\]
By a straightforward calculation we get
\[\lev u(R,n)\,|\,\mu\des \ = \ \sum_{\mskd} \:u_{\mskd}\,\lev \psi_{\mskd}\,|\,
 \partial f\des \ = \  h^d\, \sum_{\mskd} \, u_{\mskd}\,\mu_{\mskd}
 \ = \ h^d\, \lev {\bf u} \,|\, \msmu_n\des_R,\]
proving Assertion. {\QED}

In other words the sequence of functions $\mu(n) = \Phi_n(R) \msmu_n(R) \in E_n(
R,{\bbR}^d)$ converges strongly in $W_2^{-1}({\bbR}^d)$ to $\mu$. This property
implies the convergence of numbers $\lev v(n)|\mu(n)\des$ to the number
$\lev v|\mu\des$ for any $W_2^1$-weakly convergent sequence of functions
$v(n)$.

From the expression in proof of this Lemma we can get the estimate
\[ |\lev{\bf u}_n|\msmu_n\des_R| =  h^{-d}|(\partial u(R,n)|f)|
 \leq \ h^{-d} \Vert \mu \Vert_{2,-1} 
 \Vert\partial u(R,n)\Vert_2 \leq h^{-d/2}\,\Vert \mu \Vert_{2,-1}\, 
 q_R( {\bf u}_n)^{1/2},\]
where the last inequality follows from Theorem~\ref{thn3.1}. The obtained
inequality is valid for any $R$, implying:
\begin{equation}\label{exp5.5}
 \lnorm \msmu_n \lnorm_{avg,2,-1} \ = \ \frac{1}{vol(R)}\,\sum_R\,
 \lnorm \msmu_n \lnorm_{R2,-1} \ \leq \ h^{-d/2}\, \Vert \mu  \Vert_{2,-1}.
\end{equation}

Inequality (\ref{exp5.5}) and the variational equalities (\ref{exp5.1}) imply the first
result towards our proof of convergence of approximate solutions. If
${\bf u}_n$ solve (\ref{exp5.1}) then
\begin{equation}
 \Md\:\lnorm {\bf u}_n \lnorm_{avg,2,1}^2 \ \leq \
 \lev {\bf u}_n\,|\,(\lambda I+A_n)\,{\bf u}_n \des \ \leq \  
 \lnorm {\bf u}_n \lnorm_{avg,2,1}\,\lnorm \msmu_n \lnorm_{avg,2,-1}.
\end{equation}

\begin{corollary}\label{corn5.2} Let ${\bf u}_n = R(\lambda,A_n)\msmu_n$ and
$u(R,n) = \Phi_n(R) {\bf u}_n$. Then
for each $R$ the sequence $\mathfrak{U} = \{ u(R,n): n \in {\bbN}\} 
\subset \cup_n E_n(R,{\bbR}^d)$ converges weakly in $W_2^1({\bbR}^d)$ to some
$u \in W_2^1({\bbR}^d)$.
\end{corollary}

Let $u^\ast$ be the solution of (\ref{ex2.5}). Then the sequence of functions 
$\hat{u}^\ast(n)$, defined by~(\ref{exn3.3}), strongly converges to $u^\ast$ in $W_2^1$. In the 
remaining part of this analysis we have to demonstrate the expected property 
$\lim_n u(R,n) = \lim_n \hat{u}^\ast(R,n) = u^\ast$ for each $R$. We follow the 
well-known finite element technique.

\begin{equation}\begin{array}{rcl}
 \Md\, h^d\,\lnorm {\bf u}_n-\hat{{\bf u}}_n^\ast \lnorm_{avg,2,1}^2 &\leq&
 h^d\,\lev{\bf u}_n-\hat{{\bf u}}_n^\ast \,|\,(\lambda I+A_n)\,
 ({\bf u}_n-\hat{{\bf u}}_n^\ast )\des \\  &=&
 h^d\,\lev {\bf u}_n-\hat{{\bf u}}_n^\ast\, |\,(\lambda I+A_n)\, {\bf u}_n\des \\  &-&
 h^d\,\lev {\bf u}_n-\hat{{\bf u}}_n^\ast \,|\,(\lambda I+A_n)\,\hat{{\bf u}}_n^\ast \des \\
 = h^d\,\ \lev {\bf u}_n-\hat{{\bf u}}_n^\ast\,|\, \msmu_n\des  &-&
 h^d\,\lev {\bf u}_n-\hat{{\bf u}}_n^\ast \,|\,(\lambda I+A_n)\,\hat{{\bf u}}_n^\ast \des.
 \end{array}
\end{equation}
By Lemma \ref{lem5.2} the first term on the right hand side converges to 
$\lev u-u^\ast |\mu\des$. By the consistency property of Proposition 
\ref{Prop5.1} the second term converges to the same value. 

\begin{theorem}\label{Th5.1} Let $\mathfrak{U}$ be as in Corollary \ref{corn5.2}.
Then the sequence $\mathfrak{U}$ converges $W_2^1({\bbR})$-strongly to the unique 
solution $u^\ast$ of (\ref{ex2.5}).
\end{theorem}

From this result, Lemma \ref{lem5.2} and Lemma \ref{lem3.1} we get another important 
result for $\lambda = 0$.

\begin{corollary}\label{corn5.3} 
Let $D$ be a bounded domain with Lipsithz boundary and $\mu \in \dot{W}_2^{-1}(D)$.
Let $A_n(D)$ be the restriction to $G_n(D)$ of $A_n$, $\msmu_n$ on
$G_n(D)$ satisfy (\ref{exp5.3}) and ${\bf u}_n = A_n(D)^{-1}\msmu_n$.
Then the sequence $\mathfrak{U}$ converges 
strongly in $W_2^1({\bbR}^d)$ to the unique weak solution $u$ of (\ref{ex2.5}).
\end{corollary}

\subsection{Convergence in $C$- and $L_1$-spaces}
Regularity properties of solutions of the converging sequences $\mathfrak{U}$
of Corollary \ref{corn5.3} follow from the well-known
results of DeGiorgi type. We can use criteria developed by \cite{LU} as in our
approach \cite{LR3} in order to get the following result. If $\Vert \mu \Vert_\infty
=1$ then there exists $\alpha > 0$ such that $\mathfrak{U}$ of Corollary \ref{corn5.3}
converges in the H\"{o}lder spaces $C^{(\alpha)}(\overline{D})$ to $u$ of (\ref{ex2.5}).
Finally, if $\mu \in L_1(D)$, or $\mu \in {\cal P}(D)$,
numerical solutions must converge in the Banach space
$\dot{L}^{(\alpha)}(D)$, which is defined in~\cite{St}. This space is the
completion of $\dot{C}^{(\alpha)}(\overline{D})$ in the norm:
\[ \Vert \,u\,\Vert \ = \ \Vert \,u\,\Vert_1 \:+\:
 \sup_{|\mbox{\boldmath $w$}| > 0} \: \frac{1}{|\mbox{\boldmath $w$}|^\alpha}
 \Vert \, \big(Z(\mbox{\boldmath $w$}) - I\big) \,u\,\Vert_1. \]
Due to the fact that the spaces $E_n(D,{\bbR}^d)$ are finite-dimensional,
we can use the duality of finite-dimensional spaces $l_1(G_n(D))$ and 
$l_\infty(G_n(D))$ and get a result on the convergence in $L_1(D)$-space.
In order to describe the convergence in H\"{o}lder space and $L_1(D)$-space
we have to explain some details about discretizations. For the discussion of
convergence
in H\"{o}lder space we assume $\mu = f \in L_\infty(D)$, $\Vert f \Vert_\infty=1$,
while for the convergence in $L_1(D)$-space we assume $\mu \in {\cal P}(D)$.
Discretizations of $f$ are denoted by ${\bf f}_n$ and it is supposed that the
sequence of functions $f(n) = \Phi {\bf f}_n$ converges strongly to $f$ in
$L_2(D)$. The grid solutions ${\bf u}_n = A_n^{-1}{\bf f}_n$ define a
sequence of approximate solutions $\mathfrak{U} = \{ u(n): n \in {\bbN}\} \subset 
\cup_n E_n(D)$. Discretizations of $\mu \in {\cal P}(D)$ are denoted by $\msmu_n$
and constructed so that $\mu(n) \in L_1(D)$, $\lim_n \Vert \mu(n)\Vert_1 = 1$,
as well as $\lim_n \lev g|\mu(n)\des = \lev g|\mu\des$ for each $g \in
\dot{C}(D)$. The corresponding grid-solutions ${\bf v}_n = A_n^{-1}\msmu_n$ 
define a sequence of approximate solutions $\mathfrak{V} = \{ v(n): n \in {\bbN}\} 
\subset \cup_n E_n(D)$. 
\begin{corollary}\label{cor5.4} 
The sequence $\mathfrak{U}$ has a subsequence converging to $u$ in
$\dot{C}^{(\alpha)}(\overline{D})$ for some $\alpha > 0$. The sequence 
$\mathfrak{V}$ has a subsequence converging to $u$ in $L_1^{(\beta)}(D)$ 
for some $\beta > 0$.
\end{corollary}
\section{Numerical realizations}\label{sec6}
The convergence proofs of previous section were given for a simple case of only
one $\msr \in {\bbN}^d$ realizing the compartmental structure of system matrices.
In applications the coefficients $a_{ij}, i \ne j$ vary and we can gain the
compartmental structure by using several parameters $\msr$. The set ${\bbR}^d$ as
well as $D \subset {\bbR}^d$ must be decomposed into subsets $D_l$ which are
defined during description of basic schemes. It is assumed that the compartmental
structure of system matrices $A_n$ can be attained by a single parameter $\msr(l)$
for each $G_n(l) = G_n \cap \overline{D_l}$. Because of a numerous 
technical details, we cannot go into a thorough analysis of the construction of
system matrices for variable $\msr$. Nevertheless, in the following example
this approach is undertaken, and seems to be natural. 

\begin{example}
\end{example}
In this example we consider the differential operator $A = -\sum_{ij = 1}^2 \partial_i
a_{ij}\partial_j$ with the diffusion tensor defined as follows
\[ a \ = \ \left [ \begin{array}{cc} \sigma^2 & \alpha(\msx) \\ \alpha(\msx) & 1
 \end{array} \right ],  \quad \alpha(\msx) \:=\:\rho {\bbJ}_{D_0}(\msx), \quad
 \rho^2 < \sigma^2,\]
where $\sigma^2$ is a positive number, $\rho$ is a real number and $D_0 =
(1/4,3/4)^2$. 

Let $D = (0,1)^2 \subset {\bbR}^2$ and $\partial D$ be its boundary.
The function $\msx \mapsto u^\ast(x_1,x_2) = x_1x_2$ is the unique solution 
to the following boundary value problem
\[\begin{array}{lll}
 \big(A u \big)(\msx)\ = \ f(\msx) &{\rm for}& \msx \in D,\\
 u|\partial D \ = \  u^\ast|\partial D, &&\end{array} \]
where
\[ f(\msx) \ = \  2\,\rho \,{\bbJ}_{D_0}(\msx) \,+\,\frac{\rho}{4}
 \Big[ \delta \Big(x_1-\frac{1}{4}\Big)- 3\delta\Big(
 x_1-\frac{3}{4}\Big)+ \delta\Big(x_2- \frac{1}{4}\Big)-3 \delta\Big(x_2-
 \frac{3}{4}\Big)\Big].\]
The set ${\bbR}^2$ is discretized by the grid $G_n$ of grid-knots $\msx_{kl} = 
hk\mse_1 + hl\mse_2, k,l \in {\bbZ}$, where $h$ is a grid-step. In order to get
a discretization of $D$ suitable for numerical handling we assume $h =1/N$
where $N = 4M$. In this way we have discretizations $G_n(D)$ of the open
square $D = (0,1)^2$ defined by grid-knots $\msx_{kl} = (x_k, y_l), 1 \leq k,l < N$. 
The sets $G_n(2) = \overline{D_0}\cap G_n$ and $G_n(1) = G_n \setminus 
G_n(0)$ define a partition of $G_n$.
The set $G_n(2)$ is closed, while $G_n(1)$ is open. Let $G_n(1,D) = G_n(1) \cap D$.
Then $G_n(2), G_n(1,D)$ is a partition of $G_n(D)$ to be used in constructions
of numerical grids and approximate solutions.
\begin{Figure}\label{fgn6.1}
\begin{displaymath}
	\epsfxsize=6cm
	\epsfbox{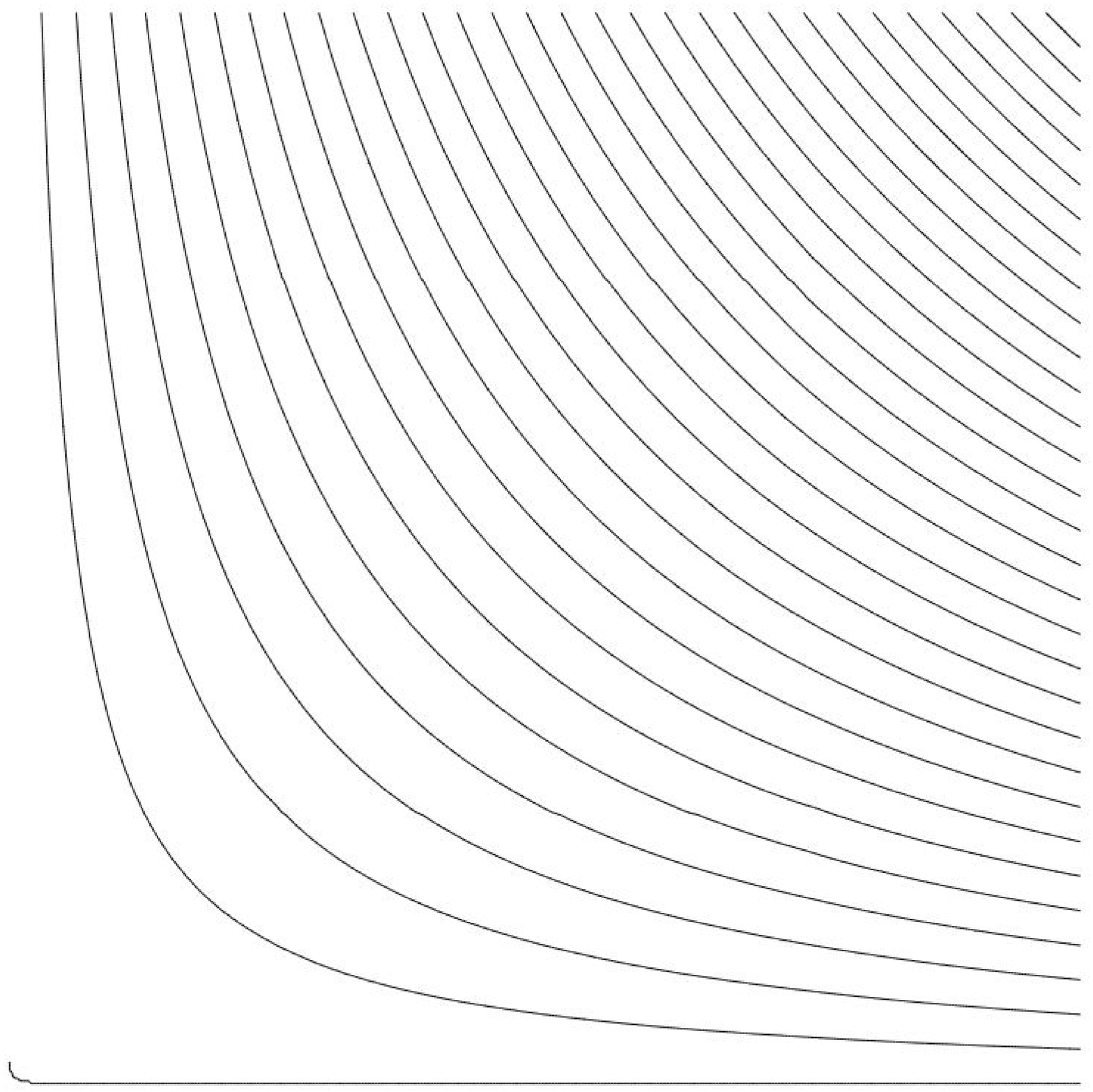}
\end{displaymath}
\end{Figure}

To demonstrate the efficiency of the extended schemes we choose 
$\sigma^2 = 10, \rho =2$ and the scheme parameters 
$r_1 = 1, r_2 =3$ as illustrated in Figure~\ref{fgn4.1}, part (c).
These values of parameters ensure the compartmental structure of 
the system matrix.

For a numerical illustration we discretize the domain $D$  by $399 \times 399$
grid-knots with the grid-step $h = 1/400$, and define the 
grid-solution by the linear system~(\ref{ex2.10}). 
This system is then solved by iterations as follows.
Let $K_n = {\rm diag}(A_n)$ and $Q_n= A_n - {\rm diag}(A_n)$. 
Then $Q_n \geq 0$ and
\begin{equation}\label{exs3.2}
 A_n^{-1} \ = \ \sum_{k=0}^\infty \: \big ( K_n^{-1} Q_n \big)^k\,
 K_n^{-1}. 
\end{equation}
Let ${\bf u}_n(r)$ be the approximation of ${\bf u}_n$  after
$r$ iterations. By taking the stopping criteria to be 
$\lnorm {\bf u}_n(r+1) - {\bf u}_n(r) \lnorm_1 < 10^{-9}$, 
we have found that the iteration terminates after $r=210$ iterations.
Then we compared the numerical approximations ${\bf u}_n(210)$ and values of the solution 
in $L_\infty(D)$ and $L_1(D)$--norms. The $L_\infty(D)$--norm of the
difference $u(h\msk) - ({\bf u}_n(210))_{\mskd}$ was $0.004$, realized at the grid-knot 
with coordinates $(0.75,0.75)$ which is one of the discontinuity points for 
$a_{12}$. The relative error in $L_1(D)$--norm was estimated by
\[ \varepsilon_{rel} \ = \ 100\:\frac{\lnorm {\bf u}^\ast \,-\,{\bf u}_n(210)\lnorm_1}{
 \lnorm {\bf u}^\ast \lnorm_1},\]
where $\big({\bf u}^\ast\big)_{\mskd} = u(h\msk)$. We obtained 
$\varepsilon_{rel} = 0.0008$, and the corresponding approximate solution $u(n)$
as illustrated in the Figure~\ref{fgn6.1}.

\end{document}